\documentclass[11pt]{article}
\usepackage[a4paper,margin=1in]{geometry}
\usepackage{amsmath,amssymb,amsthm,mathtools,bm}
\usepackage{aliascnt}
\usepackage{microtype}
\usepackage{enumitem}
\usepackage{needspace}
\usepackage{booktabs}
\usepackage{graphicx}
\usepackage{placeins}
\usepackage{tikz}
\usepackage{comment}
\usetikzlibrary{arrows.meta,calc,positioning,decorations.pathreplacing}
\usepackage{float}
\usepackage{algorithm}
\usepackage{algpseudocode}
\usepackage[authoryear,round]{natbib}
\setcitestyle{aysep={,},citesep={,}}
\usepackage[colorlinks=true,citecolor=blue,linkcolor=black,urlcolor=black]{hyperref}
\usepackage[nameinlink,capitalize,noabbrev]{cleveref}
\usepackage{lmodern}
\usepackage[T1]{fontenc}
\usepackage{setspace}
\allowdisplaybreaks

\newtheorem{theorem}{Theorem}[section]
\newaliascnt{lemma}{theorem}
\newtheorem{lemma}[lemma]{Lemma}
\aliascntresetthe{lemma}
\newaliascnt{proposition}{theorem}

\aliascntresetthe{proposition}
\newaliascnt{corollary}{theorem}
\newtheorem{corollary}[corollary]{Corollary}
\aliascntresetthe{corollary}
\theoremstyle{definition}
\newaliascnt{assumption}{theorem}
\newtheorem{assumption}[assumption]{Assumption}
\aliascntresetthe{assumption}
\newaliascnt{definition}{theorem}

\aliascntresetthe{definition}
\theoremstyle{remark}
\newaliascnt{remark}{theorem}
\newtheorem{remark}[remark]{Remark}
\aliascntresetthe{remark}

\crefname{theorem}{Theorem}{Theorems}
\Crefname{theorem}{Theorem}{Theorems}
\crefname{lemma}{Lemma}{Lemmas}
\Crefname{lemma}{Lemma}{Lemmas}
\crefname{proposition}{Proposition}{Propositions}
\Crefname{proposition}{Proposition}{Propositions}
\crefname{corollary}{Corollary}{Corollaries}
\Crefname{corollary}{Corollary}{Corollaries}
\crefname{assumption}{Assumption}{Assumptions}
\Crefname{assumption}{Assumption}{Assumptions}
\crefname{definition}{Definition}{Definitions}
\Crefname{definition}{Definition}{Definitions}
\crefname{remark}{Remark}{Remarks}
\Crefname{remark}{Remark}{Remarks}

\newcommand{\E}{\mathbb E}
\newcommand{\Pp}{\mathbb P}
\newcommand{\R}{\mathbb R}
\newcommand{\1}{\mathbf 1}
\newcommand{\Leb}{\operatorname{Leb}}
\newcommand{\dist}{\operatorname{dist}}
\newcommand{\TV}{\operatorname{TV}}
\newcommand{\KL}{\operatorname{KL}}

\newcommand{\sgn}{\operatorname{sgn}}
\newcommand{\cQ}{\mathcal Q}
\newcommand{\cP}{\mathcal P}
\newcommand{\cE}{\mathcal E}
\newcommand{\cD}{\mathcal D}

\newcommand{\wh}{\widehat}
\newcommand{\argmin}{\operatorname*{arg\,min}}

\title{Nonparametric Classification under Adversarial Attacks: Decision-Boundary Geometry and Optimal Rates}
\author{
Yanhang Zhang$^{1}$\thanks{These authors contributed equally to this work.},
Jingfu Peng$^{1}$\footnotemark[1],
and Yuhong Yang$^{1,2}$\\[0.8em]
{\small $^{1}$Yau Mathematical Sciences Center, Tsinghua University.}\\
{\small $^{2}$Beijing Institute of Mathematical Sciences and Applications.}
}
\date{}

\begin{document}
\maketitle

\begin{abstract}
Unlike in continuous regression, a small input perturbation
in classification can move an input across a decision boundary
and flip the predicted label. Decision-boundary geometry is therefore
essential to understanding adversarial classification risk.
Under anisotropic smoothness and margin conditions, we
quantify this geometric contribution and establish sharp
minimax rates for adversarial excess risk relative to the
ordinary Bayes risk.
These rates comprise a statistical term governed by the harmonic
mean smoothness and a geometric cost of adversarial perturbations
governed by the least smooth direction. As smoothness and the margin exponent vary, the
adversarial term exhibits a phase transition from a fractional
power of the attack radius to linear scaling, with a logarithmic correction in a critical regime.
In addition, we construct an adaptive classifier that attains
these rates over a range of smoothness and margin classes.
We further characterize the minimax rate for adversarial excess
risk under a budget on ordinary excess risk, yielding an explicit
trade-off curve. The curve identifies attack regimes where ordinary
minimax optimality entails a suboptimal adversarial rate.
Finally, numerical experiments support the theoretical findings.
\end{abstract}

{\small\noindent\textbf{Keywords.}
Adversarial classification; anisotropic smoothness; minimax rates; decision-boundary geometry; accuracy--robustness trade-off.\par}

\section{Introduction}
Adversarial attacks pose a serious threat to the reliability of
classification systems, particularly in applications where prediction
errors can compromise safety or security.
Early studies showed that nearly imperceptible input perturbations
can cause highly accurate image classifiers to make incorrect
predictions \citep{Szegedy2014,Goodfellow2015}.
Subsequent work has demonstrated adversarial vulnerabilities in
traffic-sign recognition \citep{EykholtEtAl2018},
speaker recognition
\citep{ChenEtAl2023QFA2SR},
and text sentiment classification
\citep{XuEtAl2024PromptAttack}.
This vulnerability has motivated the development of adversarial
training and defenses
\citep[e.g.,][]{Madry2018,ZhangEtAl2019}.

Much of the statistical theory for adversarial classification has
been developed under parametric models. Their distributional
assumptions seem to be restrictive and may not be flexible
enough to reflect the nature of data in modern applications.
This motivates nonparametric approaches to adversarial learning.
Recent work has established minimax rates for nonparametric
regression under adversarial attacks, providing insight into how
smoothness, sample size, and attack strength affect estimation
accuracy. However, in our context,
classification differs fundamentally from regression.
Even a correctly classified input can incur adversarial loss
if the classifier predicts the opposite label within the
attack radius.
Decision-boundary geometry determines which inputs are
susceptible to such label reversals.
We therefore seek to bound the probability mass of susceptible
inputs by analyzing this geometry and combine the result with
bounds on regression estimation error to control adversarial
excess risk.

We address the following three questions.
\begin{itemize}
    \item How does decision-boundary geometry affect adversarial classification risk?

    \item What are the optimal rates under anisotropic smoothness and margin conditions, and can they be attained adaptively?

    \item Is there a fundamental trade-off between ordinary and adversarial risks?
\end{itemize}
Experiments in image classification show that sensitivity to
adversarial perturbations can vary substantially across different
directions \citep{YinEtAl2019Fourier,TsuzukuSato2019Fourier}.
Motivated by this directional variation, we study these questions
in an anisotropic nonparametric setting, allowing the regression
function to have different degrees of H\"older smoothness across
coordinates.

\subsection{Boundary susceptibility probability}
\label{sec:bsp}

Let $(X,Y)$ have distribution $P$ on $[0,1]^d\times\{0,1\}$,
with covariate marginal $P_X$ and regression function
$\eta(x):=P(Y=1\mid X=x)$.
For a radius $r\ge0$, define the attack set
$A_r(x):=\{x'\in[0,1]^d:\|x'-x\|_2\le r\}$.
The worst-case adversarial risk of a classifier $f$ is defined as
\citep{Madry2018,ZhangEtAl2019}
\[
L_r(f):=\E\left\{\sup_{x'\in A_r(X)}\1\{f(x')\ne Y\}\right\}.
\]
At $r=0$, this is the ordinary misclassification risk
$L_0(f)=P\{f(X)\ne Y\}$, minimized by the Bayes classifier
$f^*(x):=\1\{\eta(x)\ge1/2\}$.
We measure adversarial excess risk relative to the ordinary Bayes
risk, that is, $\cE_r(f):=L_r(f)-L_0(f^*)$.
This benchmark retains the \(r\)-dependent excess of the optimal
adversarial risk over the ordinary Bayes risk even when the
distribution is known, as shown in \cref{cor:population}.

Given a classifier $f$, the boundary susceptibility probability (BSP) is
\[
B_r(f):=P_X\{\exists x'\in A_r(X):f(x')\ne f(X)\}.
\]
It is the probability that $f$ takes both labels on $A_r(X)$,
so an admissible perturbation can cause misclassification
for either response label.
For fixed $f$ and $r$, BSP depends only on the decision regions
of $f$ and $P_X$.
\Cref{thm:decomposition} decomposes adversarial excess risk,
up to universal constants, into ordinary excess risk and BSP as
\[
\cE_r(f)\asymp \cE_0(f)+B_r(f).
\]
The main challenge is to obtain sharp bounds on $B_r(f)$
under anisotropic smoothness and margin conditions.

\subsection{Related work}

A substantial line of work develops statistical theory for adversarial
classification under parametric distributional models.
\citet{JavanmardSoltanolkotabi2022} provide a precise asymptotic
analysis of adversarially trained linear classifiers under
high-dimensional Gaussian mixture models, characterizing standard and
robust classification accuracies in a proportional-growth regime.
\citet{DanEtAl2020} establish minimax guarantees for adversarial
excess risk relative to the optimal adversarial risk under Gaussian
mixture models.
Related analyses study whether standard and robust risks
can be minimized simultaneously.
\citet{TsiprasEtAl2019} show that, under a specific data distribution,
high standard classification accuracy necessarily entails
low robust classification accuracy.
\citet{DobribanEtAl2023} derive optimal adversarially robust classifiers
in Gaussian mixture settings and show that class imbalance can prevent
the standard and robust risks from being minimized simultaneously.
\citet{JavanmardMehrabi2024} show that low-dimensional latent structure
can mitigate this trade-off, allowing classifiers to be nearly optimal
for both standard and robust accuracy.

For nonparametric regression under adversarial attacks,
\citet{PengYang2024} establish minimax rates under adversarial
sup-norm loss.
\citet{PengYang2025} extend this analysis to adversarial $L_q$ risks
for $1\le q<\infty$ and construct adaptive estimators that attain
the minimax rates up to a logarithmic factor.
More recently, \citet{PengYang2026} show that exact interpolation
can be suboptimal for adversarial learning even under small
future $X$-attacks.
For isotropic H\"older smoothness $\beta$, our local polynomial
classifier uses a bandwidth with an $r^{1/(\beta+1)}$ term,
compared with the $r$ term in the regression analysis of
\citet{PengYang2025}.
This larger scale provides stronger smoothing to control
boundary susceptibility.

Existing work has established robustness and consistency guarantees for nonparametric classification under adversarial attacks.
\citet{WangJhaChaudhuri2018} analyze the robustness of
$k$-nearest-neighbor classifiers and show how it depends on the
choice of $k$.
They also construct a modified $1$-nearest-neighbor rule with
asymptotic robustness guarantees.
\citet{YangEtAl2020} develop generic attacks and defenses for
nonparametric classifiers and introduce adversarial pruning as a
finite-sample approximation to an optimally robust classifier.
\citet{BhattacharjeeChaudhuri2020} show that nearest-neighbor
and kernel classifiers can achieve robust consistency
when combined with adversarial pruning.

Several studies provide minimax rates and nonasymptotic risk
bounds for adversarial classification with different attack
models and risk benchmarks.
For isotropic H\"older smoothness $0<\beta\le1$
and margin exponent $\alpha$, \citet{XingEtAl2021} establish
the minimax rate
$n^{-\beta(1+\alpha)/(2\beta+d)}\vee r^{\beta(1+\alpha)}$
under random test perturbations and a prescribed attack based
on the true $\eta$, both specified independently of the fitted
classifier $f$.
Their excess risk relative to the ordinary Bayes risk equals
$\E\{|2\eta(X)-1|\,\1\{f(X')\ne f^*(X)\}\}$, where $X'$ denotes the perturbed input.
The minimax analysis uses the margin condition to handle both
estimation and perturbation.
In contrast, our worst-case analysis requires a separate sharp
geometric bound on $B_r(f)$ in addition to control of ordinary
excess risk.
Under worst-case attacks, \citet{Dohmatob2024} derives bounds
for linear classifiers that are minimax optimal up to
logarithmic factors, using the best robust linear risk as the
benchmark rather than the ordinary Bayes risk.
Relatedly, \citet{liu2024nonasymptotic} derive nonasymptotic
upper bounds for adversarial excess risk relative to a robust
oracle, accounting for the error in approximating
H\"older-smooth target functions by neural networks.
Their adversarial classification bound requires
$|\eta(x)-1/2|>c$ for all $x$ and some fixed $c>0$.
Our analysis allows $\eta$ to cross $1/2$ and does not restrict
the decision boundary to a hyperplane.

Classical nonparametric classification theory relates
minimax rates to the complexity of decision sets
\citep{MammenTsybakov1999}.
In adversarial classification,
\citet{GarciaTrillosMurray2022} derive geometric evolution
equations for smooth optimal decision boundaries, while
\citet{BungertEtAl2023} establish a nonlocal perimeter
formulation and the existence of minimizers with regular
boundaries.
Under suitable regularity conditions,
\citet{MorrisMurray2026} establish Hausdorff convergence of
optimal adversarial decision regions to the Bayes region
on compact sets as the attack radius vanishes.
\citet{LangeVasilyan2025} bound adversarial risk by ordinary
risk plus boundary volume, which corresponds to
$B_r(f)$ defined here.
Under subgaussian isotropic log-concave marginals, their
learning algorithm simultaneously guarantees, with high
probability, ordinary risk at most $O(\varepsilon)$ above
that of the best linear classifier and boundary volume
$O(r+\varepsilon)$, where $\varepsilon$ is the target accuracy.
Our analysis complements these studies by quantifying how
decision-boundary geometry affects minimax rates under
anisotropic H\"older smoothness and margin conditions.

\subsection{Organization}
The remainder of the paper is organized as follows.
Section~\ref{sec:setup} gives the loss decomposition and specifies
the model assumptions. Section~\ref{sec:minimax} develops the
minimax theory through matching upper and lower bounds. Section~\ref{sec:adaptive} establishes a holdout oracle inequality
and applies it to adaptive bandwidth selection. Section~\ref{sec:budget} characterizes the minimax
trade-off under an ordinary-risk budget. Section~\ref{sec:numerical}
presents numerical experiments, and
Section \ref{sec:conclusion} concludes with a discussion. The proofs and auxiliary
technical results are collected in the appendices.

\section{Problem formulation}
\label{sec:setup}
Throughout this paper, all classifiers are Borel measurable maps from $[0,1]^d$ to $\{0,1\}$.

\subsection{A geometric characterization of adversarial excess risk}
Recall that
\[
\cE_0(f)
=2\E_X\Bigl\{|\eta(X)-\tfrac12|\1\{f(X)\ne f^*(X)\}\Bigr\}.
\]
The following decomposition holds for every classifier $f$ and distribution $P$.
\begin{theorem}
\label{thm:decomposition}
For any classifier $f$ and $r\ge0$,
\begin{equation}
\frac12\{\cE_0(f)+B_r(f)\}
\le
\cE_r(f)
\le
\cE_0(f)+B_r(f).
\label{eq:decomposition}
\end{equation}
\end{theorem}
\begin{remark}
\label{rem:decomposition-comparison}
\citet{ZhangEtAl2019} and \citet{BungertEtAl2023} give exact
risk decompositions which yield
\[
\cE_r(f)=\cE_0(f)+
P\Bigl\{
f(X)=Y,\ \exists x'\in A_r(X):
f(x')\ne f(X)
\Bigr\}.
\]
This boundary term counts only correctly classified observations
whose predictions can be changed by an attack.
By comparison, $B_r(f)$ includes susceptibility at both correctly and
incorrectly classified observations and, for fixed $f$ and $r$, depends
only on $P_X$.
The two-sided comparison in Theorem~\ref{thm:decomposition}
allows $B_r(f)$ to be used directly
in both the minimax upper and lower bounds.
\end{remark}

\subsection{Model class and assumptions}
For $\bm\beta\in(0,\infty)^d$, write
\[
\bar\beta:=\frac{d}{\sum_{j=1}^d\beta_j^{-1}},\qquad
\beta_-:=\min_{1\le j\le d}\beta_j,\qquad
\beta_+:=\max_{1\le j\le d}\beta_j.
\]
Let $\mathbb N_0:=\{0,1,\ldots\}$.
For $\bm k=(k_1,\ldots,k_d)\in\mathbb N_0^d$, let
$
|\bm k|:=\sum_{j=1}^d k_j.
$
Here $e_j$ is the $j$th standard basis vector of $\mathbb R^d$,
and $\partial_j=\partial/\partial x_j$.
We assume that $\eta\in\mathcal C^{\bm\beta}(L)$, where
$\mathcal C^{\bm\beta}(L)$ is the coordinatewise anisotropic
H\"older class specified in
\citet{ShenGaoWittenHan2020}.

\begin{assumption}[Anisotropic H\"older regularity]
\label{ass:holder}
For each $j$, let
$
m_j:=\lceil\beta_j\rceil-1.
$
Let $\mathcal C^{\bm\beta}(L)$ be the class of continuous functions $\eta$ on $[0,1]^d$ such that every section in
coordinate $j$
has $m_j$ continuous derivatives and satisfies
\[
\|\partial_j^\ell \eta\|_{\infty,[0,1]^d}\le L,
\qquad \ell=1,\ldots,m_j,
\]
and
\[
\left|
\partial_j^{m_j}\eta(x+te_j)-\partial_j^{m_j}\eta(x)
\right|
\le L|t|^{\beta_j-m_j},
\qquad x,x+te_j\in[0,1]^d.
\]
\end{assumption}

\begin{remark}
\label{rem:anisotropic-isotropic}
When $\beta_1=\cdots=\beta_d=\beta$, \cref{ass:holder} is equivalent,
up to constants, to the usual isotropic H\"older condition
for noninteger $\beta$ and for $\beta=1$
\citep{AudibertTsybakov2007}.
At integer orders $\beta\ge2$, the coordinatewise condition
is in general weaker than the corresponding isotropic condition.
\end{remark}

\begin{assumption}[Margin condition]
\label{ass:margin}
For some $\alpha>0$ and $C_\alpha>0$,
\[
P_X\left\{\left|\eta(X)-\frac12\right|\le t\right\}
\le C_\alpha t^\alpha,
\qquad t>0.
\]
\end{assumption}

\begin{assumption}[Strong density]
\label{ass:density}
$P_X$ has support $[0,1]^d$ and a density $\mu$ satisfying
\[
0<\mu_-\le\mu(x)\le\mu_+<\infty,
\qquad x\in[0,1]^d.
\]
\end{assumption}

\begin{assumption}[Parameter range]
\label{ass:range}
The parameters satisfy
\[
\alpha\bar\beta\le d,
\qquad
\alpha(1\wedge\beta_+)\le1.
\]
\end{assumption}
\begin{remark}
The condition $\alpha(1\wedge\beta_+)\le1$ is used in our
multivariate lower bound, which perturbs $\eta$ in every
coordinate direction.
In the isotropic case,
\citet{AudibertTsybakov2007} show that
$\alpha(1\wedge\beta)>1$ prevents $\eta$ from crossing $1/2$
in the interior of the covariate support.
With our support and density assumptions, a constant Bayes
classifier then attains the ordinary Bayes risk for every
attack radius.
\end{remark}

Let $\cP(\alpha,\bm\beta,L)$ denote the class of distributions
satisfying \cref{ass:holder,ass:margin,ass:density,ass:range},
with fixed constants $C_\alpha$, $\mu_-$ and $\mu_+$.
For a sample $\mathcal D_n=((X_1,Y_1),\ldots,(X_n,Y_n))\sim P^{\otimes n}$,
define
\[
\mathcal R_{n,r}
:=
\inf_{\wh f}
\sup_{P\in\cP(\alpha,\bm\beta,L)}
\E_{P^{\otimes n}}\cE_r(\wh f).
\]

\section{Minimax rates}
\label{sec:minimax}
We establish matching upper and lower bounds for $\mathcal R_{n,r}$.
The upper bound is obtained using a piecewise local polynomial
classifier. The lower bound follows from separate constructions
for the statistical and adversarial terms.
\subsection{Piecewise local polynomial classifier and upper bound}
\label{sec:estimator}
Fix an integer $p\ge1$ with $\beta_+<p+1$.
For $\bm h=(h_1,\ldots,h_d)$ satisfying $h_j^{-1}\in\mathbb N$
for all $j$, let $\cQ_{\bm h}$ be the regular partition of $[0,1]^d$ with side lengths $\bm h$.
For each $Q\in\cQ_{\bm h}$, let $z_Q$ be its center.
Let $\underline h:=\min_{1\le j\le d}h_j$ be the minimum bandwidth,
$V_{\bm h}:=\prod_{j=1}^d h_j$ be the cell volume, and
$b_{\bm h}:=\sum_{j=1}^d h_j^{\beta_j}$ be the approximation scale.
For $x\in Q$, define the basis vector by
\[
\phi_Q(x)
:=
\left(
\prod_{j=1}^d
\left(\frac{x_j-z_{Q,j}}{h_j}\right)^{k_j}
\right)_{\bm k\in\mathbb N_0^d,\ |\bm k|\le p}.
\]
The empirical Gram matrix and response vector are
\[
\wh G_Q
:=
\frac1{nV_{\bm h}}
\sum_{i=1}^n
\1\{X_i\in Q\}\phi_Q(X_i)\phi_Q(X_i)^\top,
\qquad
\wh b_Q
:=
\frac1{nV_{\bm h}}
\sum_{i=1}^n
\1\{X_i\in Q\}\phi_Q(X_i)Y_i.
\]
Define
\[
\kappa_p:=\mu_-\lambda_{\min}\left(
\left[
\int_{[-1/2,1/2]^d}\prod_{j=1}^d u_j^{k_j+k'_j}\,du
\right]_{|\bm k|,|\bm k'|\le p}
\right)>0,
\]
and
\[
\widehat\theta_Q
:=
\begin{cases}
\widehat G_Q^{-1}\widehat b_Q,
&
\lambda_{\min}(\widehat G_Q)\ge \kappa_p/2,
\\[1mm]
(1/2,0,\ldots,0)^\top,
&
\text{otherwise}.
\end{cases}
\]
The fitted polynomial is
\[
\wh P_Q(x):=\wh\theta_Q^\top\phi_Q(x),
\qquad x\in Q.
\]
For $x\in Q$, the local polynomial classifier is
\[
\wh f_{\bm h}(x):=\1\{\wh P_Q(x)\ge1/2\}.
\]
For $0<t\le1$, define
\begin{equation}
\Lambda_{\alpha,\beta_-}(t)
:=
\begin{cases}
t^{\alpha\beta_--1},&\alpha\beta_-<1,\\[1mm]
\log(e/t),&\alpha\beta_-=1,\ \beta_->1,\\[1mm]
1,&\text{otherwise}.
\end{cases}
\label{eq:Lambda}
\end{equation}

\Needspace{14\baselineskip}
\begin{theorem}
\label{thm:fixed-upper}
Under \cref{ass:holder,ass:margin,ass:density,ass:range}, there exists
$C>0$ such that, for $nV_{\bm h}\ge1$ and
$0\le r\le\underline h$,
\begin{align}
\sup_{P\in\cP(\alpha,\bm\beta,L)}
\E\cE_r(\wh f_{\bm h})
\le C\Bigg\{&
b_{\bm h}^{1+\alpha}
+(nV_{\bm h})^{-(1+\alpha)/2}
\notag\\*
&+r\Lambda_{\alpha,\beta_-}(\underline h)
\left\{1+\frac{b_{\bm h}+(nV_{\bm h})^{-1/2}}
{\underline h^{\beta_-}}\right\}^{1+\alpha}
+\1\{\beta_->1\}\left(\frac r{\underline h}\right)^2
\Bigg\}.
\label{eq:fixed-upper}
\end{align}
\end{theorem}
In \eqref{eq:fixed-upper}, the terms $b_{\bm h}^{1+\alpha}$
and $(nV_{\bm h})^{-(1+\alpha)/2}$ bound the approximation
and stochastic estimation contributions to ordinary excess risk,
respectively.
The last two terms bound $\E B_r(\wh f_{\bm h})$.
When $\alpha\beta_-<1$, the factor
$r\underline h^{\alpha\beta_--1}$ in the boundary bound increases
as $\underline h$ decreases.
Balancing this contribution with the approximation term
$b_{\bm h}^{1+\alpha}$ gives the $r^{1/(\beta_-+1)}$ term
in the minimum bandwidth of \cref{cor:oracle-upper}.

For $0<r\le1$, define
\begin{equation}
\psi_{\alpha,\beta_-}(r):=
\begin{cases}
r^{\frac{\beta_-(1+\alpha)}{\beta_-+1}},
& \alpha\beta_-<1,\\[0.6em]
r\log(e/r),
& \alpha\beta_-=1,\ \beta_->1,\\[0.6em]
r,
& \alpha\beta_->1
\ \text{or}\ \alpha\beta_-=1,\ \beta_-\le1.
\end{cases}
\label{eq:psi-definition}
\end{equation}
Optimizing \eqref{eq:fixed-upper} over the bandwidths gives the
following upper bound, with $\psi_{\alpha,\beta_-}(r)$ as the
adversarial term.

\begin{corollary}
\label{cor:oracle-upper}
Under \cref{ass:holder,ass:margin,ass:density,ass:range}, there exist bandwidths satisfying
\[
h_j\asymp
n^{-\frac{\bar\beta}{\beta_j(2\bar\beta+d)}}
\vee r^{\frac{\beta_-}{\beta_j(\beta_-+1)}},
\qquad j=1,\ldots,d,
\]
such that
\begin{equation}
\sup_{P\in\cP(\alpha,\bm\beta,L)}
\E\cE_r(\wh f_{\bm h})
\le
C\left\{
n^{-\frac{\bar\beta(1+\alpha)}{2\bar\beta+d}}
+
\psi_{\alpha,\beta_-}(r)
\right\}.
\label{eq:oracle-upper}
\end{equation}
\end{corollary}
In the isotropic case, the bandwidth choice in
\cref{cor:oracle-upper} reduces to
\[
h\asymp n^{-1/(2\beta+d)}\vee r^{1/(\beta+1)}.
\]
The bandwidth term involving $r$ is of order
$r^{1/(\beta+1)}$, compared with order $r$ in adversarial
nonparametric regression \citep{PengYang2025}. In classification, the need to control boundary susceptibility probability leads
to a larger bandwidth scale and hence stronger smoothing than
in regression.

\subsection{Geometric mechanism for the lower bound}
\begin{figure}[H]
\centering

\definecolor{eraseBlue}{HTML}{075BFF}
\definecolor{eraseRed}{HTML}{FF2020}
\definecolor{eraseGreen}{HTML}{16813A}
\definecolor{eraseBlueTint}{HTML}{E5EFFF}
\definecolor{eraseGreenTint}{HTML}{E3F7E8}
\definecolor{eraseRedTint}{HTML}{FDE7E9}
\definecolor{eraseBlueEdge}{HTML}{67A0F3}
\definecolor{eraseGreenEdge}{HTML}{62B77A}
\definecolor{eraseRedEdge}{HTML}{EE9199}
\begin{tikzpicture}[
    >=Latex,
    font=\small,
    x=0.82cm,
    y=4.0cm,
    scale=0.90,
    transform shape,
    axis/.style={->, line width=0.8pt},
    reg/.style={eraseBlue, line width=1.25pt},
    clf/.style={eraseRed, dashed, line width=1.15pt},
    threshold/.style={gray!75, dashed, line width=0.8pt},
    bandblue/.style={fill=eraseBlueTint,draw=eraseBlueEdge,dashed,line width=0.5pt},
    bandgreen/.style={fill=eraseGreenTint,draw=eraseGreenEdge,dashed,line width=0.5pt},
    bandred/.style={fill=eraseRedTint,draw=eraseRedEdge,dashed,line width=0.5pt},
]

\begin{scope}[shift={(0,0)}]

\node[font=\large] at (4.15,1.22)
    {\textbf{(a) Erase:} $f$ constant};

\draw[axis] (0,-0.08) -- (8.55,-0.08) node[right] {$t$};
\draw[axis] (0,-0.08) -- (0,1.15);

\node[left] at (0,0) {$0$};
\node[left] at (0,0.5) {$1/2$};
\node[left] at (0,1) {$1$};

\def\xLp{1.05}
\def\xLpend{2.25}
\def\xLm{6.00}
\def\xLmend{7.20}

\fill[bandblue]  (\xLp,0) rectangle (\xLpend,1);
\fill[bandgreen] (\xLm,0) rectangle (\xLmend,1);

\draw[threshold] (0,0) -- (8.4,0);
\draw[threshold] (0,0.5) -- (8.4,0.5);

\draw[reg,domain=0:8.3,samples=120,smooth,variable=\x]
    plot ({\x},{0.01 + 0.99/(1 + exp(1.20*(\x-4.28)))});

\node[eraseBlue] at (4.90,0.69) {$\eta(t)$};
\draw[->,eraseBlue,line width=0.7pt]
    (4.63,0.66) -- (4.38,0.53);

\draw[clf] (0,1) -- (8.3,1);
\node[eraseRed] at (7.35,1.07) {$f(t)=1$};

\node[eraseBlue, font=\large] at ({(\xLp+\xLpend)/2},0.30) {$J^+$};
\node[eraseGreen, font=\large] at ({(\xLm+\xLmend)/2},0.64) {$J^-$};

\draw[<->,line width=0.85pt]
    (1.65,0.5) -- (1.65,0.98)
    node[midway,right] {$h^{\beta_-}$};

\draw[<->,line width=0.85pt]
    (6.50,0.5) -- (6.50,0.10)
    node[midway,right] {$h^{\beta_-}$};

\draw[line width=0.8pt]
    (\xLp,-0.115) -- (\xLp,-0.165)
    (\xLpend,-0.115) -- (\xLpend,-0.165)
    (\xLp,-0.145) -- (\xLpend,-0.145);
\node[below] at ({(\xLp+\xLpend)/2},-0.15) {$h$};

\draw[line width=0.8pt]
    (\xLm,-0.115) -- (\xLm,-0.165)
    (\xLmend,-0.115) -- (\xLmend,-0.165)
    (\xLm,-0.145) -- (\xLmend,-0.145);
\node[below] at ({(\xLm+\xLmend)/2},-0.15) {$h$};

\end{scope}

\begin{scope}[shift={(10.05,0)}]

\node[font=\large] at (4.15,1.22)
    {\textbf{(b) Retain:} $f$ nonconstant};

\draw[axis] (0,-0.08) -- (8.55,-0.08) node[right] {$t$};
\draw[axis] (0,-0.08) -- (0,1.15);

\node[left] at (0,0) {$0$};
\node[left] at (0,0.5) {$1/2$};
\node[left] at (0,1) {$1$};

\def\xRp{0.75}
\def\xRpend{1.95}
\def\xRm{6.00}
\def\xRmend{7.20}

\pgfmathsetmacro{\xb}{4.18 + ln(50/49)/1.20}
\pgfmathsetmacro{\xAttackL}{\xb-0.425}
\pgfmathsetmacro{\xAttackR}{\xb+0.425}

\fill[bandblue]  (\xRp,0) rectangle (\xRpend,1);
\fill[bandgreen] (\xRm,0) rectangle (\xRmend,1);
\fill[bandred]   (\xAttackL,0) rectangle (\xAttackR,1);

\draw[black!65,dashed,line width=0.8pt] (\xb,-0.08) -- (\xb,1);

\draw[threshold] (0,0.5) -- (8.4,0.5);

\draw[reg,domain=0:8.3,samples=120,smooth,variable=\x]
    plot ({\x},{0.01 + 0.99/(1 + exp(1.20*(\x-4.18)))});

\node[eraseBlue] at (5.05,0.67) {$\eta(t)$};
\draw[->,eraseBlue,line width=0.7pt]
    (4.78,0.64) -- (4.52,0.51);

\draw[clf] (0,1) -- (\xb,1);
\draw[clf] (\xb,0) -- (8.3,0);

\node[eraseRed] at (1.95,1.07) {$f(t)=1$};
\node[eraseRed] at (7.55,0.15) {$f(t)=0$};

\draw[line width=0.8pt]
    (\xAttackL,1.015) -- (\xAttackL,1.065)
    (\xAttackR,1.015) -- (\xAttackR,1.065)
    (\xAttackL,1.04) -- (\xAttackR,1.04);
\node[above] at ({(\xAttackL+\xAttackR)/2},1.04) {$r$};

\node[eraseRed,align=center,font=\footnotesize] (vulnerable) at (2.9,0.22)
    {vulnerable\\region};

\node[eraseBlue, font=\large] at ({(\xRp+\xRpend)/2},0.30) {$J^+$};
\node[eraseGreen, font=\large] at ({(\xRm+\xRmend)/2},0.64) {$J^-$};

\draw[<->,line width=0.85pt]
    (1.45,0.5) -- (1.45,0.98)
    node[midway,right] {$h^{\beta_-}$};

\draw[<->,line width=0.85pt]
    (6.50,0.5) -- (6.50,0.10)
    node[midway,right] {$h^{\beta_-}$};

\draw[line width=0.8pt]
    (\xRp,-0.115) -- (\xRp,-0.165)
    (\xRpend,-0.115) -- (\xRpend,-0.165)
    (\xRp,-0.145) -- (\xRpend,-0.145);
\node[below] at ({(\xRp+\xRpend)/2},-0.15) {$h$};

\draw[line width=0.8pt]
    (\xRm,-0.115) -- (\xRm,-0.165)
    (\xRmend,-0.115) -- (\xRmend,-0.165)
    (\xRm,-0.145) -- (\xRmend,-0.145);
\node[below] at ({(\xRm+\xRmend)/2},-0.15) {$h$};

\end{scope}

\end{tikzpicture}

\caption{
\textbf{Erase-or-retain mechanism.}
In panel (a), the classifier is constant on $J$ and disagrees
with the Bayes rule on one of the two regions, contributing
at least $c h^{\beta_-+1}$ to ordinary excess risk.
In panel (b), the classifier takes both labels on $J$,
and the set of susceptible points in $J$ has length at least $cr$.
}
\label{fig:erase-retain}
\end{figure}
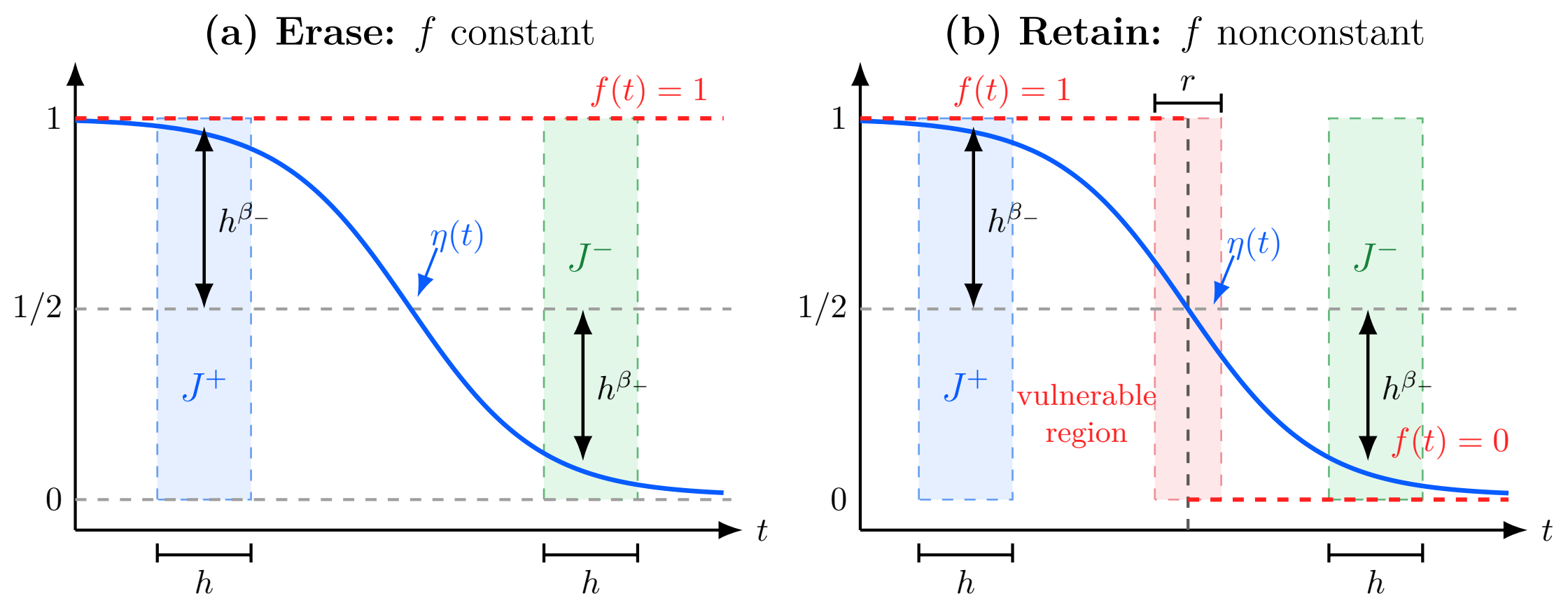

The adversarial lower bound follows from a dichotomy on
one-dimensional coordinate sections.
On an interval where $\eta$ lies on both sides of $1/2$,
a constant classifier incurs ordinary excess risk, whereas
a classifier taking both labels contributes to the BSP.
The following lemma quantifies these two contributions.
\begin{lemma}
\label{lem:erase-retain}
Suppose that $P_X$ is uniform on $[0,1]^d$ and that $\eta$
depends only on $t=x_1$. Let $J_1,\ldots,J_M\subset[0,1]$ be intervals with pairwise disjoint
interiors and $\Leb(J_k)\ge r$.
For a fixed $c_0>0$ and scales $h_1,\ldots,h_M\in(0,1]$, suppose
that each $J_k$ contains disjoint subintervals
$J_k^+$ and $J_k^-$ satisfying
\[
\Leb(J_k^+)\wedge\Leb(J_k^-)\ge c_0h_k,
\qquad
\eta\ge\frac12+c_0h_k^{\beta_-} \ \text{on }J_k^+,
\qquad
\eta\le\frac12-c_0h_k^{\beta_-} \ \text{on }J_k^-.
\]
Then, every measurable classifier $f$ satisfies
\[
\cE_0(f)+B_r(f)
\ge
c\sum_{k=1}^M\min\{h_k^{\beta_-+1},r\},
\]
where $c>0$ depends only on $c_0$.
\end{lemma}
Figure~\ref{fig:erase-retain} illustrates the two alternatives
in \cref{lem:erase-retain}.


\begin{remark}
For $\alpha\beta_-<1$, our construction uses
$M\asymp h^{\alpha\beta_--1}$ intervals at a common scale
$h\asymp r^{1/(\beta_-+1)}$.
The two costs in \cref{lem:erase-retain} then balance,
yielding a lower bound of order
$Mr\asymp r^{\beta_-(1+\alpha)/(\beta_-+1)}$.

On the critical line $\alpha\beta_-=1$, a single interval
gives an order-$r$ lower bound when $\beta_-\le1$.
For $\beta_->1$, summing order-$r$ contributions over
$\asymp\log(e/r)$ scales yields $r\log(e/r)$.
When $\alpha\beta_->1$, a single interval at a fixed scale
again gives an order-$r$ lower bound.
The constructions are given in \cref{app:lower-adversarial}.
\end{remark}

The following theorem combines the statistical and adversarial
lower bounds.
\begin{theorem}
\label{thm:lower}
Assume that \cref{ass:holder,ass:margin,ass:density,ass:range}
hold and that $C_\alpha$ is sufficiently large.
There exists $c>0$ such that, for all $n\ge1$ and $0\le r\le1$,
\begin{equation}
\mathcal R_{n,r}
\ge
c\left\{
n^{-\frac{\bar\beta(1+\alpha)}{2\bar\beta+d}}
+
\psi_{\alpha,\beta_-}(r)
\right\}.
\label{eq:lower}
\end{equation}
\end{theorem}

\begin{corollary}
\label{cor:minimax}
Under the conditions of \cref{cor:oracle-upper,thm:lower},
\[
\mathcal R_{n,r}
\asymp
n^{-\frac{\bar\beta(1+\alpha)}{2\bar\beta+d}}
+
\psi_{\alpha,\beta_-}(r).
\]
\end{corollary}

The statistical term is governed by the harmonic mean smoothness
$\bar\beta$ and reduces to the classical minimax classification rate
in the isotropic case \citep{AudibertTsybakov2007}.
The adversarial term depends on the least smooth direction
$\beta_-$ and the margin exponent $\alpha$.
Greater smoothness in the remaining coordinates can therefore
improve the statistical rate without changing the order of
the adversarial term.
The ordinary Bayes and optimal adversarial benchmarks are related by
\[
\cE_r(f)
=\left\{L_r(f)-\inf_g L_r(g)\right\}
+\left\{\inf_g L_r(g)-L_0(f^*)\right\}.
\]
The first term is excess risk relative to the optimal adversarial risk.
Our benchmark also retains the second term, the irreducible cost of
adversarial perturbations. The following corollary gives its worst-case order.

\begin{corollary}
\label{cor:population}
Under the conditions of \cref{cor:minimax}, for $0<r\le1$,
\[
\sup_{P\in\cP(\alpha,\bm\beta,L)}\inf_f\cE_r(f)
\asymp\psi_{\alpha,\beta_-}(r).
\]
\end{corollary}
The adversarial term therefore remains unavoidable in the worst case
even when the distribution is known.

\begin{figure}[H]
\centering
\definecolor{phaseBlue}{HTML}{123EDB}
\definecolor{phaseRed}{HTML}{C73713}
\definecolor{phaseRose}{HTML}{A7193D}
\definecolor{phaseGold}{HTML}{B88600}
\definecolor{phaseStatTint}{HTML}{E5F1FE}
\definecolor{phasePowerTint}{HTML}{FFF1DE}
\definecolor{phaseLinearTint}{HTML}{FDEDF2}
\definecolor{phaseExcludedTint}{HTML}{F0F0F0}
\resizebox{0.93\linewidth}{!}{%
\begin{tikzpicture}[
    x=1cm,y=1cm,>=Latex,
    every node/.style={font=\footnotesize},
    title/.style={font=\normalsize},
    region/.style={align=center,inner sep=2pt},
    balance/.style={phaseBlue,line width=1.15pt},
    critsolid/.style={phaseGold,line width=1.0pt},
    critdash/.style={black!55,dashed,line width=.85pt},
    statfill/.style={fill=phaseStatTint,draw=none},
    powerfill/.style={fill=phasePowerTint,draw=none},
    linearfill/.style={fill=phaseLinearTint,draw=none},
    exfill/.style={fill=phaseExcludedTint,draw=none}
]
\def\W{7.25}
\def\H{5.80}
\def\xcrit{4.18}
\def\ypeak{2.55}
\def\ycritical{1.96}
\def\leftbalancepath{%
    (0,1.15)
    .. controls (1.65,1.25) and (3.52,2.18) .. (\xcrit,\ypeak)
}
\def\rightbalancecontinuation{%
    .. controls (4.82,2.18) and (6.15,1.15) .. (\W,1.08)
}

\begin{scope}
  \node[title] at (.5*\W,\H+.40) {(a) $0<\alpha<1$};
  \fill[powerfill] (0,0) rectangle (\xcrit,\H);
  \fill[linearfill] (\xcrit,0) rectangle (\W,\H);
  \fill[statfill]
    \leftbalancepath \rightbalancecontinuation -- (\W,0) -- (0,0) -- cycle;
  \draw[balance] \leftbalancepath \rightbalancecontinuation;
  \draw[critdash] (\xcrit,0) -- (\xcrit,\ypeak);
  \draw[->] (0,0) -- (\W+.16,0) node[right] {$\beta_-$};
  \draw[->] (0,0) -- (0,\H+.16) node[above] {$r$};
  \node[below=4pt] at (\xcrit,0) {$1/\alpha$};

  \node[region,text=phaseRed] at (1.55,3.35)
    {$\alpha\beta_-<1$\\[4pt]
     $\mathcal R_{n,r}\asymp r^{\frac{\beta_-(1+\alpha)}{\beta_-+1}}$};
  \node[region,text=phaseRose] at (6.40,3.35)
    {$\alpha\beta_->1$\\[4pt]
     $\mathcal R_{n,r}\asymp r$};
  \node[region,text=phaseGold] (phaseCriticalLeft)
    at (\xcrit,4.35)
    {\textbf{Critical rate}\\[3pt]
     $\mathcal R_{n,r}\asymp r\log(e/r)$};
  \draw[critsolid] (\xcrit,\ycritical) -- ([yshift=-3pt]phaseCriticalLeft.south);
  \draw[critsolid] ([yshift=3pt]phaseCriticalLeft.north) -- (\xcrit,\H);
  \fill[white] (\xcrit,\ypeak) circle (2.8pt);
  \draw[black,line width=.65pt] (\xcrit,\ypeak) circle (2.8pt);
  \fill[black] (\xcrit,\ycritical) circle (2.1pt);
  \node[region,text=phaseBlue] at (2.68,.66)
    {\textbf{Statistical regime}\\[4pt]
     $\mathcal R_{n,r}\asymp n^{-\frac{\bar\beta(1+\alpha)}{2\bar\beta+d}}$};
\end{scope}

\begin{scope}[shift={(8.15,0)}]
  \node[title] at (.5*\W,\H+.40) {(b) $\alpha>1$};
  \fill[powerfill] (0,0) rectangle (\xcrit,\H);
  \fill[exfill] (\xcrit,0) rectangle (\W,\H);
  \fill[statfill] \leftbalancepath -- (\xcrit,0) -- (0,0) -- cycle;
  \draw[balance] \leftbalancepath;
  \draw[critdash] (\xcrit,0) -- (\xcrit,\ypeak);
  \draw[->] (0,0) -- (\W+.16,0) node[right] {$\beta_-$};
  \draw[->] (0,0) -- (0,\H+.16) node[above] {$r$};
  \node[below=4pt] at (\xcrit,0) {$1/\alpha$};

  \node[region,text=phaseRed] at (1.55,3.35)
    {$\alpha\beta_-<1$\\[4pt]
     $\mathcal R_{n,r}\asymp r^{\frac{\beta_-(1+\alpha)}{\beta_-+1}}$};
  \node[region,text=phaseGold] (phaseCriticalRight)
    at (\xcrit,4.35)
    {\textbf{Critical rate}\\[4pt]
     $\mathcal R_{n,r}\asymp r$};
  \draw[critsolid] (\xcrit,\ypeak) -- ([yshift=-3pt]phaseCriticalRight.south);
  \draw[critsolid] ([yshift=3pt]phaseCriticalRight.north) -- (\xcrit,\H);
  \fill[black] (\xcrit,\ypeak) circle (2.1pt);
  \node[region,text=phaseBlue] at (2.35,.66)
    {\textbf{Statistical regime}\\[4pt]
     $\mathcal R_{n,r}\asymp n^{-\frac{\bar\beta(1+\alpha)}{2\bar\beta+d}}$};
  \node[region,text=black!70] at (5.72,2.88)
    {Excluded by\\[3pt]\cref{ass:range}};
\end{scope}
\end{tikzpicture}%
}
\caption{
Phase transitions in the minimax rate.
The blue curves separate the statistical and adversarial regimes
away from the critical line $\beta_-=1/\alpha$.
In panel (a), the filled point marks the order of the transition
radius with logarithmic correction.
The open circle shows the corresponding order without this correction.
In panel (b), the critical case requires
$\beta_1=\cdots=\beta_d=1/\alpha$.
The filled point marks the balance between the statistical term
and the linear adversarial term.
}
\label{fig:beta-r-phase-diagram}
\end{figure}
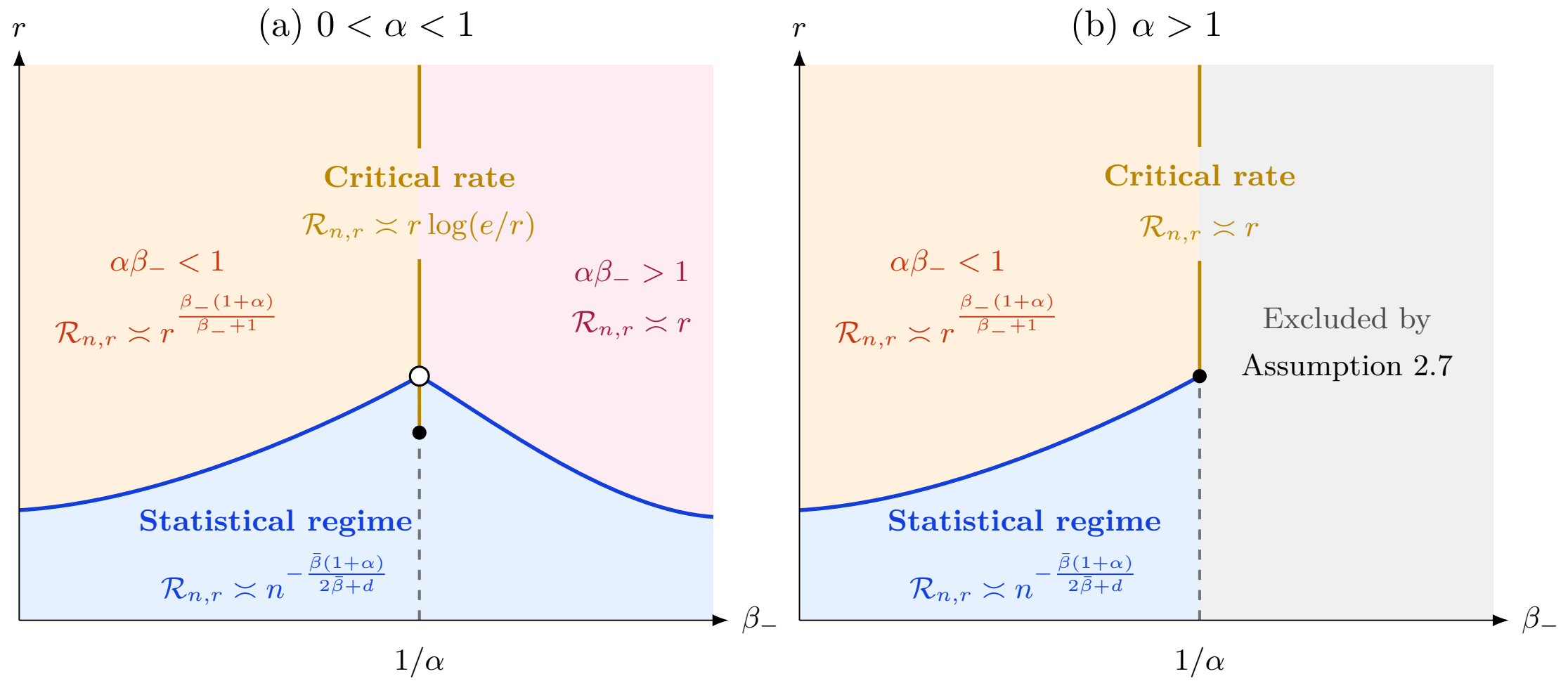

\begin{remark}[Statistical--adversarial transition]
The statistical and adversarial terms balance when
$\psi_{\alpha,\beta_-}(r)\asymp n^{-\frac{\bar\beta(1+\alpha)}{2\bar\beta+d}}$. Thus the transition radius satisfies
\[
r_n^\ast\asymp
\begin{cases}
n^{-\frac{\bar\beta(\beta_-+1)}
{\beta_-(2\bar\beta+d)}},
& \alpha\beta_-<1,\\[0.8em]
n^{-\frac{\bar\beta(1+\alpha)}
{2\bar\beta+d}}/\log n,
& \alpha\beta_-=1,\ \beta_->1,\\[0.8em]
n^{-\frac{\bar\beta(1+\alpha)}
{2\bar\beta+d}},
& \alpha\beta_->1
\ \text{or}\ 
\alpha\beta_-=1,\ \beta_-\le1.
\end{cases}
\]
For $r\ll r_n^\ast$, the usual nonparametric classification rate
\citep{AudibertTsybakov2007} dominates. For $r\gg r_n^\ast$, the adversarial term
dominates.
\end{remark}
\Cref{fig:beta-r-phase-diagram} summarizes the
transition in $r$ and the phase structure in the $(\beta_-,r)$-plane.

\section{Adaptive classification}
\label{sec:adaptive}
The rate-optimal bandwidths in \cref{cor:oracle-upper} depend on
the unknown smoothness vector. In this section, 
we construct an adaptive classifier that
attains the minimax rate without an additional logarithmic term
when $\alpha\bar\beta<d$.
\subsection{An oracle inequality for holdout selection}
\label{sec:adaptive-oracle}
Split the sample into a training sample
$
\cD_1=\bigl((X_i^{(1)},Y_i^{(1)})\bigr)_{i=1}^{n_1}
$
and an independent validation sample
$
\cD_2=\bigl((X_j^{(2)},Y_j^{(2)})\bigr)_{j=1}^{n_2},
n_1+n_2=n.
$
Let $\wh f_1,\ldots,\wh f_N$ be classifiers constructed
using $\cD_1$.

For a classifier $f$, define its empirical adversarial risk on
$\cD_2$ by
\[
\wh L_{r,2}(f)
:=
\frac1{n_2}
\sum_{j=1}^{n_2}
\1\left\{
\exists x'\in A_r(X_j^{(2)}):
f(x')\ne Y_j^{(2)}
\right\}.
\]
The holdout selector is
\[
\wh f_{r}
\in
\argmin_{f\in\{\wh f_1,\ldots,\wh f_N\}}
\wh L_{r,2}(f).
\]
For a classifier $f$, define
\[
W_{r,f}(x,y)
:=
\1\{\exists x'\in A_r(x):f(x')\ne y\}
-
\1\{f^*(x)\ne y\}.
\]

\begin{lemma}
\label{lem:bernstein}
Under \cref{ass:margin}, for every classifier $f$ and $r\ge0$,
\[
\operatorname{Var}\{W_{r,f}(X,Y)\}
\le C\cE_r(f)^{\alpha/(1+\alpha)},
\]
where $C$ depends only on $\alpha$ and $C_\alpha$.
\end{lemma}

Combining Lemma \ref{lem:bernstein} with the holdout argument
\citep{Massart2007,MaillardArlotLerasle2021} yields the following
oracle inequality.
\begin{theorem}
\label{thm:holdout}
Suppose that \cref{ass:margin} holds. For each fixed $r\ge0$ and
$\tau,\delta\in(0,1)$, conditional on $\cD_1$, with probability
at least $1-\delta$ over $\cD_2$, we have
\[
\cE_r(\wh f_{r})
\le
\frac{1+\tau}{1-\tau}
\min_{1\le m\le N}\cE_r(\wh f_m)
+
C_\tau
\left\{\frac{\log(2N/\delta)}{n_2}\right\}^{\frac{1+\alpha}{2+\alpha}}.
\]
Moreover,
\[
\E\cE_r(\wh f_{r})
\le
2\E\left\{
\min_{1\le m\le N}
\cE_r(\wh f_m)
\right\}
+
C
\left\{
\frac{\log(2N)}{n_2}
\right\}^{\frac{1+\alpha}{2+\alpha}}.
\]
The constants $C_\tau$ and $C$ depend only on
$(\tau,\alpha,C_\alpha)$ and $(\alpha,C_\alpha)$, respectively.
\end{theorem}

\subsection{Adaptive bandwidth selection}
\label{sec:vanishing-holdout}
Given an integer $p\ge1$, we consider smoothness vectors
satisfying $\beta_+<p+1$.
Take
$
n_2:=\lceil n/\log n\rceil, J_n:=\lfloor\log_2 n_1\rfloor,
$
and define
\[
\mathcal H_n
:=\left\{2^{-\bm k}:\bm k\in\mathbb N_0^d,\ |\bm k|\le J_n\right\},
\qquad
2^{-\bm k}:=(2^{-k_1},\ldots,2^{-k_d}).
\]
Every $\bm h\in\mathcal H_n$ satisfies
$n_1V_{\bm h}\ge1$. For $0\le r\le1$, let
\[
\mathcal H_n(r):=\{\bm h\in\mathcal H_n:\min_jh_j\ge r\}.
\]
For each $\bm h\in\mathcal H_n$, construct the classifier
$\wh f_{\bm h}^{(1)}$ of Section~\ref{sec:estimator}
with polynomial degree $p$ using $\cD_1$.
Apply the holdout selector to
$\{\wh f_{\bm h}^{(1)}:\bm h\in\mathcal H_n(r)\}$, namely,
\[
\wh{\bm h}(r)
\in\argmin_{\bm h\in\mathcal H_n(r)}
\wh L_{r,2}(\wh f_{\bm h}^{(1)}),
\qquad
\wh f_{\mathrm{ad},r}:=\wh f_{\wh{\bm h}(r)}^{(1)}.
\]
Combining \cref{thm:holdout,cor:oracle-upper} gives the following bound.

\begin{theorem}
\label{thm:adaptive}
Under \cref{ass:holder,ass:margin,ass:density,ass:range}, there exists
$C>0$ such that, for all $0\le r\le1$,
\[
\sup_{P\in\cP(\alpha,\bm\beta,L)}
\E\cE_r(\wh f_{\mathrm{ad},r})
\le
C\left\{
 n^{-\frac{\bar\beta(1+\alpha)}{2\bar\beta+d}}
 +\psi_{\alpha,\beta_-}(r)
 +\left(\frac{\log n\,\log\log n}{n}\right)^{\frac{1+\alpha}{2+\alpha}}
\right\}.
\]
If $\alpha\bar\beta<d$, then
\[
\sup_{P\in\cP(\alpha,\bm\beta,L)}
\E\cE_r(\wh f_{\mathrm{ad},r})
\le
C\left\{
 n^{-\frac{\bar\beta(1+\alpha)}{2\bar\beta+d}}
 +\psi_{\alpha,\beta_-}(r)
\right\}.
\]
\end{theorem}
When $\alpha\bar\beta<d$, the holdout remainder is of smaller
order than the statistical term.
Under the conditions of \cref{cor:minimax}, the selector therefore
attains the minimax rate without an additional logarithmic loss.
For comparison, the Lepski procedure of \citet{PengYang2025}
yields adaptive bounds with an additional logarithmic factor
under adversarial $L_q$ risk for $1\le q<\infty$.

\begin{remark}
For implementation, the holdout procedure can be repeated over random
splits of the same sample \citep{MaillardArlotLerasle2021}. For an odd
integer $V\ge3$, let
$\wh f_{\mathrm{ad},r}^{(1)},\ldots,\wh f_{\mathrm{ad},r}^{(V)}$
be the classifiers obtained from $V$ such splits. Define the aggregated holdout classifier by
\[
\wh f_{\mathrm{agg},r}(x)
:=\1\left\{
\frac1V\sum_{v=1}^V\wh f_{\mathrm{ad},r}^{(v)}(x)\ge\frac12
\right\}.
\]
For fixed $V$, the aggregated classifier satisfies the bounds in
\cref{thm:adaptive}, with constants depending on $V$.
\end{remark}
We refer to the resulting method as \textbf{A}dversarial \textbf{L}ocal \textbf{P}olynomial classification with \textbf{H}oldout \textbf{A}ggregation (ALPHA).
Algorithm~\ref{alg:adaptive-plp} summarizes the procedure.\begin{algorithm}[H]
\caption{ALPHA classifier}
\label{alg:adaptive-plp}
\begin{algorithmic}[1]
\Require Sample $\cD=\{(X_i,Y_i)\}_{i=1}^n$, radius $r\in[0,1]$,
degree $p\ge1$, and $V\ge1$.
\State Set $n_2=\lceil n/\log n\rceil$ and $n_1=n-n_2$.
\State Construct the admissible bandwidth set $\mathcal H_n(r)$ defined above.
\For{$v=1,\ldots,V$}
  \State Randomly split $\cD$ into $\cD_1^{(v)}$ and $\cD_2^{(v)}$ of sizes $n_1$ and $n_2$.
  \For{$\bm h\in\mathcal H_n(r)$}
    \State Fit $\wh f_{\bm h}^{(v)}$ on $\cD_1^{(v)}$ with degree $p$.
    \State Evaluate its empirical adversarial risk $\wh L_{r,2}^{(v)}(\wh f_{\bm h}^{(v)})$ on $\cD_2^{(v)}$.
  \EndFor
  \State Select $\wh{\bm h}^{(v)}(r)\in\argmin_{\bm h\in\mathcal H_n(r)}
    \wh L_{r,2}^{(v)}(\wh f_{\bm h}^{(v)})$ and set $\wh f_{\mathrm{ad},r}^{(v)}=\wh f_{\wh{\bm h}^{(v)}(r)}^{(v)}$.
\EndFor
\Ensure $\displaystyle \wh f_{\mathrm{agg},r}(x)
 =\1\!\left\{V^{-1}\sum_{v=1}^{V}\wh f_{\mathrm{ad},r}^{(v)}(x)\ge\tfrac12\right\}$.
\end{algorithmic}
\end{algorithm}

\section{Minimax trade-off between accuracy and robustness}
\label{sec:budget}
\begin{figure}[H]
\centering
\definecolor{budgetLine}{HTML}{2F78C4}
\definecolor{budgetTightTint}{HTML}{FDE8E5}
\definecolor{budgetTradeTint}{HTML}{E1F5EB}
\definecolor{budgetLooseTint}{HTML}{E8F0FF}
\resizebox{0.94\linewidth}{!}{%
\begin{tikzpicture}[
 x=1.1666667cm,y=1cm,>=Latex,
 every node/.style={font=\small},
 region/.style={align=center,inner sep=3pt},
 frontier/.style={budgetLine,line width=1.5pt},
 separator/.style={black!45,dashed,line width=.75pt}
]
\fill[budgetTightTint] (0,0) rectangle (4,5.8);
\fill[budgetTradeTint] (4,0) rectangle (8,5.8);
\fill[budgetLooseTint] (8,0) rectangle (12,5.8);
\draw[separator] (4,0) -- (4,5.8);
\draw[separator] (8,0) -- (8,5.8);
\draw[black!55,dashed,line width=.9pt]
 (0,1.196924997) -- (12,1.196924997);
\draw[frontier] (0,3.8) -- (4,3.8);
\draw[frontier,domain=4:8,samples=100,variable=\x]
 plot ({\x},{3.8*(4/\x)^(5/3)});
\draw[frontier] (8,1.196924997) -- (12,1.196924997);
\fill[budgetLine] (4,3.8) circle (2.1pt);
\fill[budgetLine] (8,1.196924997) circle (2.1pt);
\draw[->,line width=.8pt] (0,0) -- (12.3,0);
\draw[->,line width=.8pt] (0,0) -- (0,6.05);
\node[region] at (2,5.25) {\textbf{(I) Tight budget}};
\node[region] at (6,5.25)
 {\textbf{(II) Moderate budget}};
\node[region] at (10,5.25)
 {\textbf{(III) Loose budget}};
\node[region,text=budgetLine,font=\normalsize,anchor=center] at (2,3.15) {$r^{\alpha\beta_-}$};
\node[region,text=budgetLine,font=\normalsize,anchor=center] at (6,3.45)
 {$r\varepsilon^{-\frac{1-\alpha\beta_-}{\beta_-(1+\alpha)}}$};
\node[region,text=budgetLine,font=\normalsize,anchor=center] at (10,2.25)
 {$r^{\frac{\beta_-(1+\alpha)}{\beta_-+1}}$};
\draw (4,.07) -- (4,-.07);
\draw (8,.07) -- (8,-.07);
\node[below=7pt] at (4,0) {$r^{\beta_-(1+\alpha)}$};
\node[below=7pt] at (8,0)
 {$r^{\frac{\beta_-(1+\alpha)}{\beta_-+1}}$};
\node at (6,-1.32) {Budget $\varepsilon$};
\node[rotate=90] at (-.62,2.9) {Attack contribution};
\end{tikzpicture}%
}
\caption{Three budget regimes for the attack contribution in
\eqref{eq:budget-frontier}.
The horizontal dashed line indicates the order of the attack
contribution to the unconstrained minimax rate.}
\label{fig:budget-phase-transition}
\end{figure}
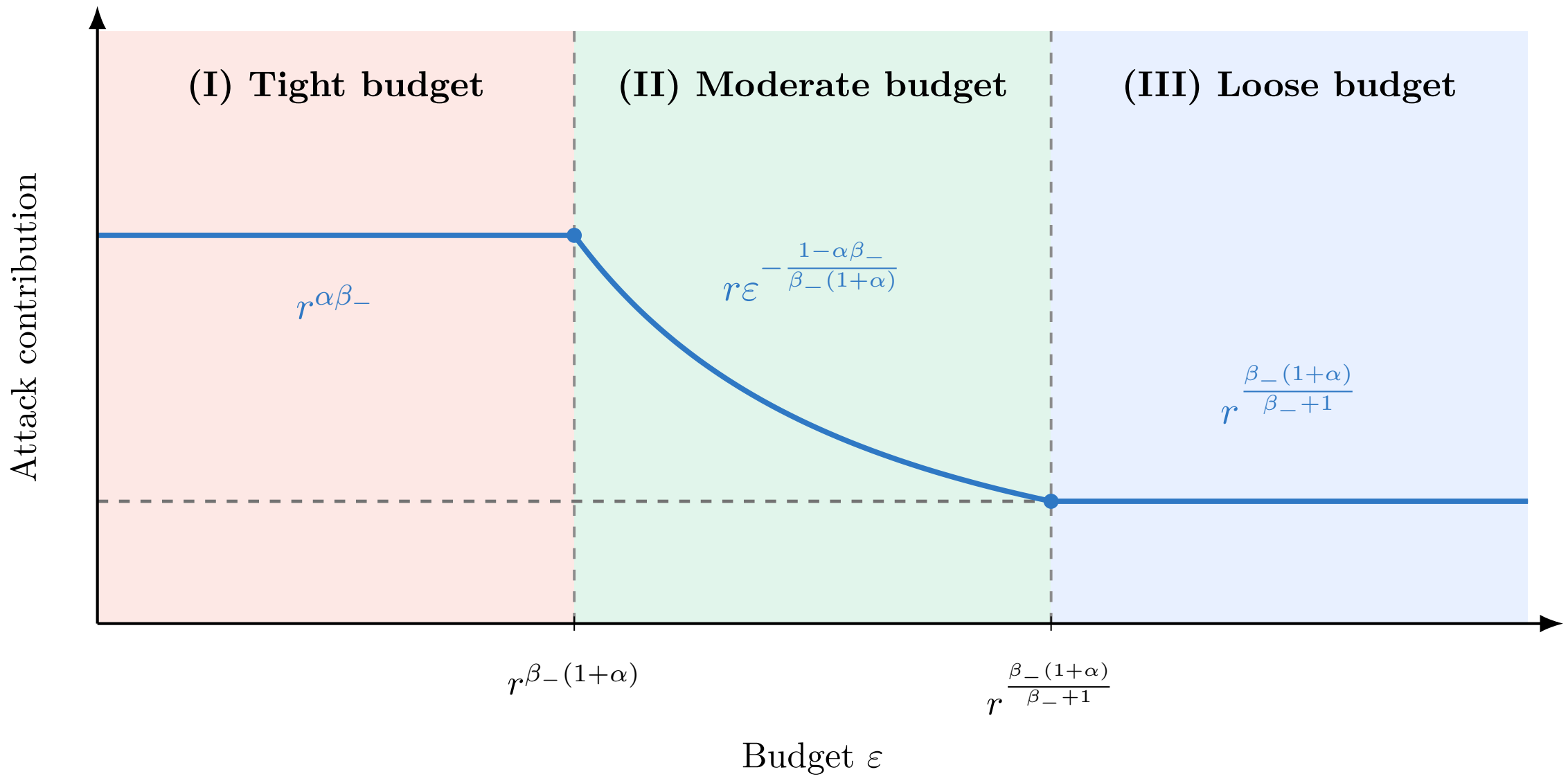

We next study the minimax adversarial excess risk under a uniform
constraint on the expected ordinary excess risk.
For a budget $\varepsilon>0$, define
\begin{equation}
\mathcal R_{n,r}(\varepsilon)
:=
\inf_{\substack{\wh f:
 \sup_{P\in\cP}\E_{P^{\otimes n}}\cE_0(\wh f)
 \le\varepsilon}}
\sup_{P\in\cP}
\E_{P^{\otimes n}}\cE_r(\wh f).
\label{eq:budget-minimax}
\end{equation}
When \(\varepsilon\) is of the ordinary minimax order, \(\mathcal R_{n,r}(\varepsilon)\) characterizes the best adversarial rate achievable while maintaining ordinary minimax optimality.
\begin{theorem}
\label{thm:budget-frontier}
Assume that \cref{ass:holder,ass:margin,ass:density,ass:range}
hold, $\beta_-\le1$, $\alpha\beta_-<1$ and $C_\alpha$ is sufficiently large.
There exists a constant $C_0>0$ such that, for all
$0\le r\le1$ and
\[
\varepsilon\ge
C_0n^{-\frac{\bar\beta(1+\alpha)}{2\bar\beta+d}},
\]
we have
\begin{equation}
\mathcal R_{n,r}(\varepsilon)\asymp
n^{-\frac{\bar\beta(1+\alpha)}{2\bar\beta+d}}
+
\begin{cases}
r^{\alpha\beta_-},
& \varepsilon\le r^{\beta_-(1+\alpha)},\\[0.9em]
r\varepsilon^{-\frac{1-\alpha\beta_-}{\beta_-(1+\alpha)}},
& r^{\beta_-(1+\alpha)}<\varepsilon\le
  r^{\frac{\beta_-(1+\alpha)}{\beta_-+1}},\\[0.9em]
r^{\frac{\beta_-(1+\alpha)}{\beta_-+1}},
& \varepsilon>r^{\frac{\beta_-(1+\alpha)}{\beta_-+1}}.
\end{cases}
\label{eq:budget-frontier}
\end{equation}
\end{theorem}
For a sufficiently small constant $c>0$, the admissible set in
\eqref{eq:budget-minimax} is empty when
$\varepsilon<
c n^{-\bar\beta(1+\alpha)/(2\bar\beta+d)}$.
The attack contribution has the following three regimes.
\begin{enumerate}[label=(\roman*),leftmargin=*]
\item \textbf{Tight budget.}
If $\varepsilon\lesssim r^{\beta_-(1+\alpha)}$, the attack
contribution has order $r^{\alpha\beta_-}$.

\item \textbf{Moderate budget.}
If
\[
r^{\beta_-(1+\alpha)}\lesssim\varepsilon
\lesssim r^{\frac{\beta_-(1+\alpha)}{\beta_-+1}},
\]
the attack contribution decreases as $\varepsilon$ increases,
giving an explicit trade-off between ordinary accuracy and
adversarial robustness.

\item \textbf{Loose budget.}
If $\varepsilon\gtrsim
r^{\beta_-(1+\alpha)/(\beta_-+1)}$, we recover the unconstrained minimax rate
\[
\mathcal R_{n,r}(\varepsilon)\asymp
n^{-\frac{\bar\beta(1+\alpha)}{2\bar\beta+d}}
+\psi_{\alpha,\beta_-}(r).
\]
\end{enumerate}
In the moderate regime, a larger budget permits larger bandwidths.
This reduces the bound on BSP but increases
the approximation term in the ordinary excess risk bound.
\Cref{fig:budget-phase-transition} illustrates the three regimes.

The following corollary identifies a range of attack radii for which
the ordinary and adversarial minimax rates cannot be attained
simultaneously.
\begin{corollary}
\label{cor:budget-conflict}
Under the conditions of \cref{thm:budget-frontier}, let
$
\varepsilon_n
=C n^{-\frac{\bar\beta(1+\alpha)}{2\bar\beta+d}},
$
where $C>0$ is any sufficiently large fixed constant.
If
\[
n^{-\frac{\bar\beta(\beta_-+1)}{\beta_-(2\bar\beta+d)}}
\ll r_n\ll1,
\]
then
\[
\frac{\mathcal R_{n,r_n}(\varepsilon_n)}{\mathcal R_{n,r_n}}
\longrightarrow\infty.
\]
\end{corollary}
Over this range of radii, a classifier attaining the ordinary
minimax rate uniformly over $\cP$ cannot also attain
the unconstrained adversarial minimax rate.
For example, let $d=1$, $\beta_-=1$, $\alpha=1/2$, and
$r=n^{-1/2}$.
For a sufficiently large constant $C$, the budget
$\varepsilon=Cn^{-1/2}$ gives
\[
\mathcal R_{n,r}\asymp n^{-3/8},
\qquad
\mathcal R_{n,r}(\varepsilon)\asymp n^{-1/3}.
\]
Both excess risks converge to zero, but the constrained risk
is larger by a factor of order $n^{1/24}$.

\section{Numerical experiments}
\label{sec:numerical}
We study bandwidth selection and the accuracy--robustness trade-off,
then compare ALPHA with other robust classifiers. In the implementation, ALPHA averages the local linear and
quadratic estimates with equal weights and combines the
classifiers from $V=11$ holdout splits by majority vote. We report averages over 100
replications at each attack radius in
$\mathcal R=\{0,0.005,\ldots,0.08\}$. Shaded bands are pointwise 95\%
Monte Carlo intervals. Paired differences are computed within each
replication as the comparator's risk minus ALPHA's risk.

\subsection{Simulation studies}
\label{sec:num-simulations}
We generate independent observations from
\[
X\sim\operatorname{Unif}([0,1]^d),
\qquad
Y\mid X=x\sim\operatorname{Ber}(\eta(x)).
\]
\subsubsection{Adversarial versus ordinary selection}
\label{sec:num-selection}
For $n=400$, consider the signed quadratic regression function
\[
\eta(x)=\frac12+\frac14(2x-1)|2x-1|.
\]
We compare classifiers selected by adversarial validation
at radius $r$ with those selected by ordinary validation
at $r=0$.

\begin{figure}[H]
  \centering
  \includegraphics[width=\textwidth]{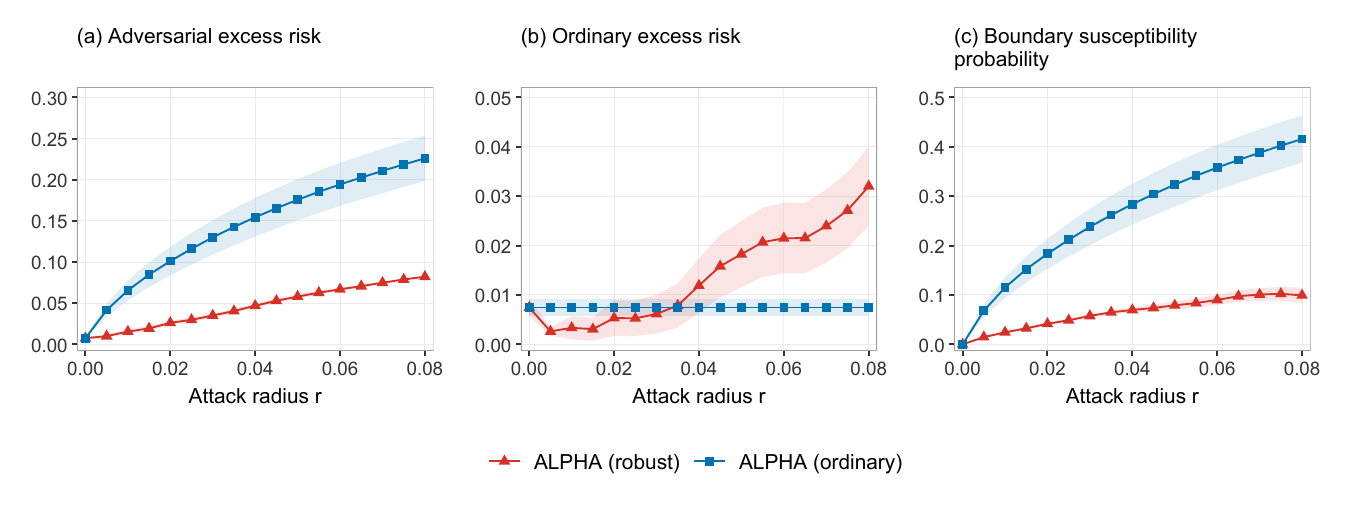}
  \caption{Adversarial excess risk, ordinary excess risk, and boundary
  susceptibility probability under adversarial and ordinary bandwidth selection.}
  \label{fig:num-selection}
\end{figure}
In Figure~\ref{fig:num-selection}, adversarial selection shows a growing advantage over ordinary selection in adversarial excess risk as the attack radius increases. Its ordinary excess risk increases at larger radii, but only modestly. In contrast, BSP rises sharply under ordinary selection and remains much higher than under adversarial selection.

\subsubsection{Anisotropic versus isotropic bandwidth selection}
\label{sec:num-anisotropy}
Let
\[
\eta(x_1,x_2)=\frac12+\frac{1+x_2}{4}
\sgn\{\sin(5\pi x_1)\}|\sin(5\pi x_1)|^{1/2}.
\]
The regression function is rough in the first coordinate and affine in
the second. For $n\in\{200,400,600\}$, we compare anisotropic bandwidth
selection with the isotropic restriction $h_1=h_2$.

\begin{figure}[H]
  \centering
  \includegraphics[width=\textwidth]{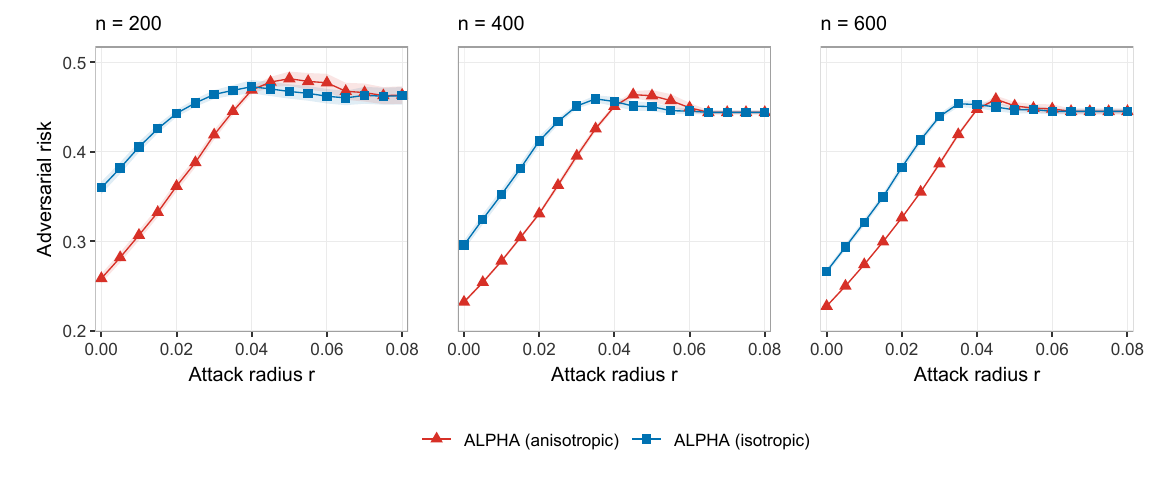}
  \caption{Adversarial risk of ALPHA with anisotropic and isotropic bandwidth selection.}
  \label{fig:num-anisotropy}
\end{figure}
Figure~\ref{fig:num-anisotropy} shows that anisotropic selection
yields noticeably lower adversarial risk at small radii across all sample sizes.
This advantage diminishes as $r$ increases, and the two
methods have similar risks at large radii.

\subsubsection{Accuracy--robustness trade-off}
\label{sec:sim-budget}
Set $n=3000$ and let
\[
\eta(x)=\frac12+\frac{x-1/4}{2}
\left\{1-2e^{-(16x-8)^2}\right\}
\left\{1-2e^{-(16x-12)^2}\right\}.
\]
In each replication, we form an oracle benchmark by minimizing the
numerically evaluated $\cE_r(f)$ over the library of fitted classifiers
subject to $\cE_0(f)\le\varepsilon$.
We consider $\varepsilon\in\{0.020,0.025,\ldots,0.100\}$ and
$r\in\{0.02,0.04,0.08\}$.

\begin{figure}[htbp]
  \centering
  \includegraphics[scale=0.8]{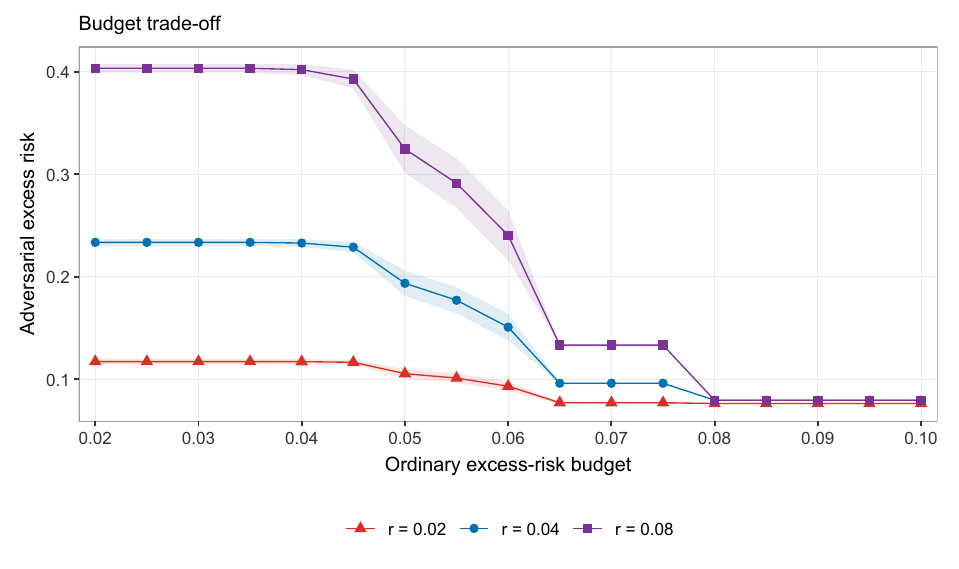}
  \caption{Adversarial excess risk under an ordinary excess-risk budget.}
  \label{fig:sim-budget}
\end{figure}
In Figure~\ref{fig:sim-budget}, adversarial excess risk remains nearly constant for small
budgets $\varepsilon$, decreases over an intermediate range,
and levels off again for large budgets.
This pattern illustrates the trade-off in
Theorem~\ref{thm:budget-frontier}.

\subsubsection{Comparison with other robust classifiers}
\label{sec:num-comparison}
For $n\in\{200,400,600\}$, consider
\[
\eta(x_1,x_2)=\frac12+\frac{1+x_2}{4}
\sgn\{\cos(3\pi x_1)\}|\cos(3\pi x_1)|^{1/2}.
\]
We compare ALPHA with two nearest-neighbor methods,
AP-$k$NN \citep{YangEtAl2020} and Robust 1-NN
\citep{WangJhaChaudhuri2018}.
We also consider two methods for adversarial training of neural
networks, PGD-AT \citep{Madry2018} and TRADES
\citep{ZhangEtAl2019}.

\begin{figure}[H]
  \centering
  \includegraphics[width=\textwidth]{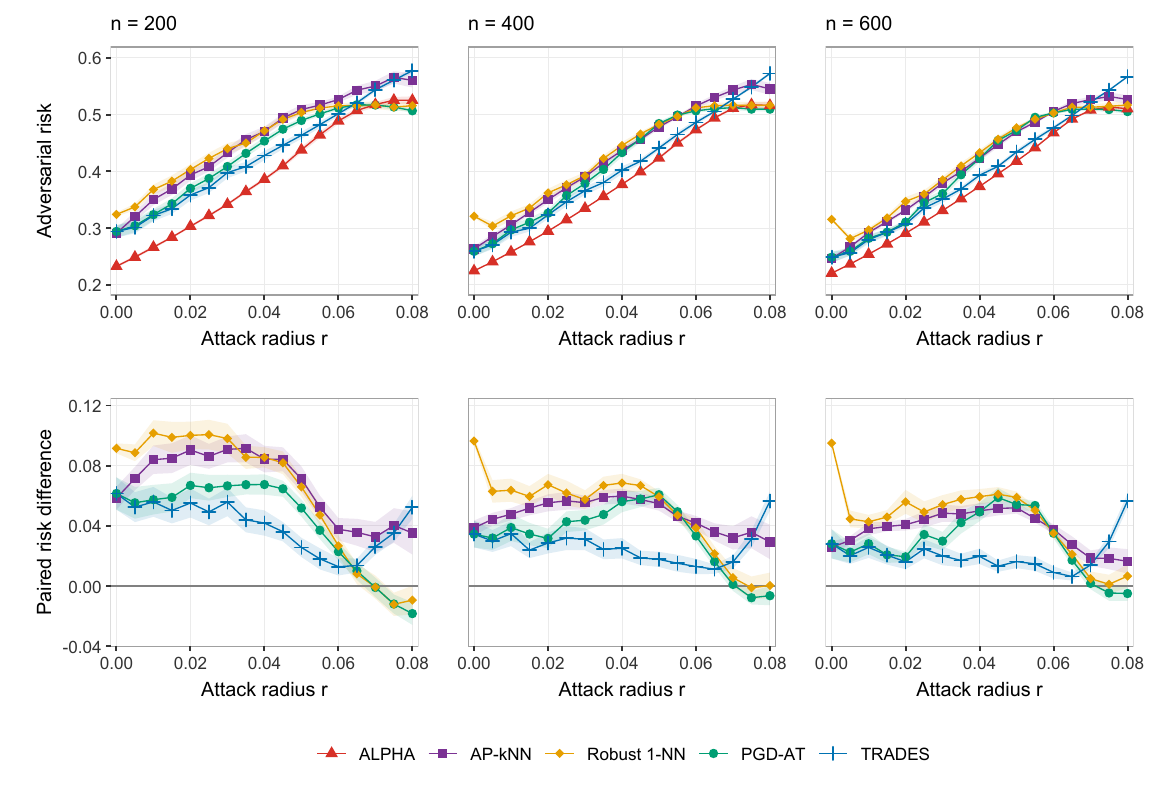}
  \caption{Adversarial risks and paired differences relative to ALPHA.}
  \label{fig:num-comparison}
\end{figure}
Figure~\ref{fig:num-comparison} shows that ALPHA has lower
adversarial risk than all four comparators at small to moderate
attack radii across all three sample sizes.
At larger attack radii, ALPHA remains competitive with the other methods.

\subsection{Breast cancer data analysis}
\label{sec:num-realdata}
The Wisconsin Diagnostic Breast Cancer (WDBC) data \citep{WDBC1993} contain 569 cases, including 357 benign and 212 malignant cases. Following \citet{MangasarianStreetWolberg1995}, we use mean texture, worst area, and worst smoothness as the three predictors.
These predictors are extracted from digitized images of cell nuclei.
We assess the robustness of benign--malignant
classification to adversarial perturbations of these features.

We randomly divide the data into an 80\% training set and a 20\%
test set using stratified sampling. Within each split, the predictors
are rank-transformed to $[0,1]$ using the training sample, and the
same transformation is applied to the test sample.
Adversarial perturbations are applied in this transformed
feature space, with the labels held fixed.
All five methods use a training attack radius of $\rho=0.08$
and are evaluated on the test sample over $\mathcal R$ using the
same attack pool. We repeat this procedure 100 times and report
the average test errors.

\begin{figure}[H]
  \centering
  \includegraphics[width=\textwidth]{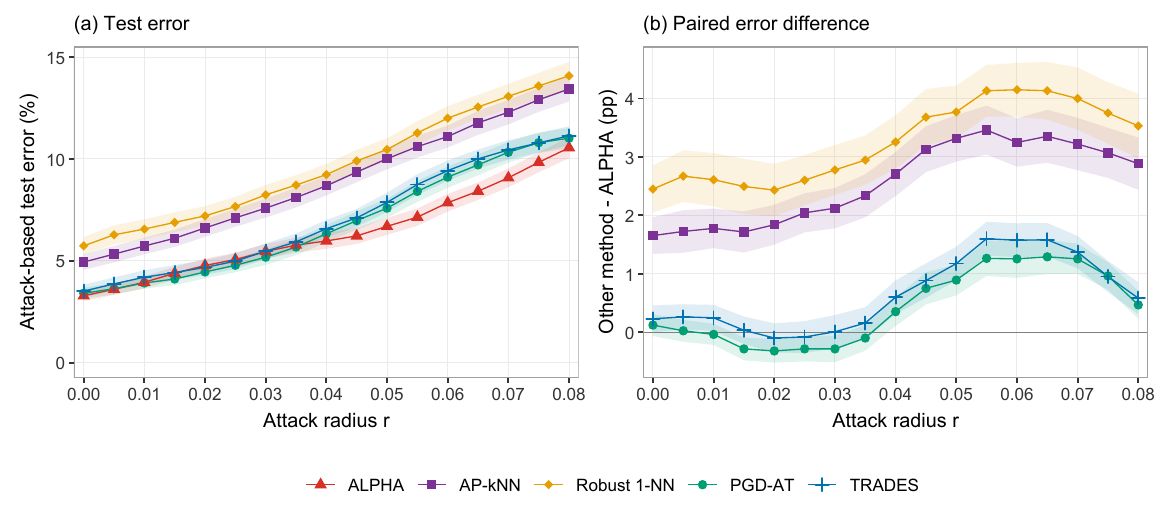}
  \caption{Test errors and paired differences relative to ALPHA on the WDBC data.}
  \label{fig:num-realdata}
\end{figure}

Figure~\ref{fig:num-realdata} shows the test-error curves and paired differences, while Table~\ref{tab:num-realdata} summarizes test errors over the full attack-radius grid.
\begin{table}[H]
  \centering
  \caption{Test errors averaged over the 17 attack radii, with training attack radius $\rho=0.08$.
  Standard errors over the 100 splits are in parentheses.}
  \label{tab:num-realdata}
  \begin{tabular}{lc}
    \toprule
    Method & Mean test error \\
    \midrule
    ALPHA & \textbf{0.0635 (0.0017)}  \\
    AP-$k$NN & 0.0892 (0.0022)  \\
    Robust 1-NN & 0.0961 (0.0024)  \\
    PGD-AT & 0.0679 (0.0018)  \\
    TRADES & 0.0701 (0.0019)  \\
    \bottomrule
  \end{tabular}
\end{table}
\Needspace{4\baselineskip}
ALPHA has the lowest average test error over all evaluated
attack radii. It is
competitive with PGD-AT and TRADES at small radii and has lower error
at larger radii. Both nearest-neighbor methods have higher errors
throughout the radius grid.
\FloatBarrier

\section{Conclusion and Discussion}
\label{sec:conclusion}

We show that adversarial excess risk is comparable to the
sum of ordinary excess risk and boundary susceptibility.
Building on this characterization, we establish sharp
minimax rates that separate a statistical term governed
by the harmonic mean smoothness from an adversarial cost
governed by the least smooth direction.
This cost remains unavoidable in the worst case even when
the distribution is known.
We also construct an adaptive classifier that attains
these rates over a range of smoothness and margin classes.

Reducing boundary susceptibility can
come at the cost of ordinary accuracy.
Our constrained minimax analysis quantifies this trade-off
through an explicit curve under a budget on ordinary
excess risk.
In certain attack regimes, ordinary minimax optimality
entails a suboptimal adversarial rate.

Future work could examine how the geometry of more general
perturbation sets affects boundary susceptibility and
minimax rates.
Another question is whether boundary susceptibility helps
explain the directional sensitivity of neural network
classifiers.

\appendix
\paragraph{Notation for the appendix.}
Let $\Delta(x):=|\eta(x)-1/2|$.
For $A\subset E$, $\partial_E A$ is the boundary of $A$ relative
to $E$, and $\partial A$ is its Euclidean boundary.
We write $\Leb$ for Lebesgue measure in the ambient dimension and
$\#A$ for the cardinality of a set $A$.
For $A\subset\mathbb R^d$ and $r\ge0$, let
\[
A^{+r}:=\{x\in\mathbb R^d:\operatorname{dist}(x,A)\le r\},
\]
where $\operatorname{dist}$ denotes Euclidean distance.
For a classifier $f$ and a fixed attack radius $r$, define
\[
V_f(x)
:=
\mathbf 1\bigl\{
\exists x'\in A_r(x): f(x')\neq f(x)
\bigr\}.
\]
\section{Proof of the loss decomposition}
\label{app:decomposition}
\begin{proof}[Proof of \cref{thm:decomposition}]
For any $x\in[0,1]^d$, since $x\in A_r(x)$, $V_f(x)=1$ if and only if
$f$ takes both values on $A_r(x)$. Define
\[
e_r(f,x)
:=
\Pp\{\exists x'\in A_r(x):f(x')\neq Y\mid X=x\}
-
\min\{\eta(x),1-\eta(x)\}.
\]
By conditioning on $X$, $\cE_r(f)=\E_X[e_r(f,X)]$.

If $V_f(x)=0$, then $f$ is constant on $A_r(x)$, so
\[
e_r(f,x)=2\Delta(x)\1\{f(x)\neq f^*(x)\}.
\]
If $V_f(x)=1$, $f$ takes both values on $A_r(x)$, so the conditional
adversarial loss equals one and
\[
e_r(f,x)=1-\left\{\frac12-\Delta(x)\right\}
=\frac12+\Delta(x).
\]
Since
$0\le2\Delta(x)\1\{f(x)\neq f^*(x)\}\le2\Delta(x)\le1$,
the two cases give
\[
\Delta(x)\1\{f(x)\neq f^*(x)\}+\frac12V_f(x)
\le
e_r(f,x)
\le
2\Delta(x)\1\{f(x)\neq f^*(x)\}+V_f(x).
\]
Integrating with respect to $P_X$ proves \eqref{eq:decomposition}.
\end{proof}

\section{Proof of the upper bound}
\label{app:upper}

\begin{figure}[H]
\centering
\begin{tikzpicture}[
 font=\small,
 box/.style={draw=black!65,rounded corners=2pt,align=center,
   inner sep=7pt,text width=.39\textwidth,minimum height=1.4cm},
 wide/.style={box,text width=.86\textwidth,fill=black!3},
 flow/.style={-{Stealth[length=2mm]},draw=black!70,rounded corners=2pt}]
 \node[wide] (start) {\textbf{Statistical--geometric decomposition}\\[3pt]
 $\E\cE_r\le\E\cE_0+\E B_r$};
 \node[box,below=7mm of start.south,xshift=-.235\textwidth] (stat)
 {\textbf{Ordinary excess risk}\\[3pt]
 Cellwise concentration + margin\\[3pt]
 $\E\cE_0\lesssim b_{\bm h}^{1+\alpha}+v_{\bm h}^{1+\alpha}$};
 \node[box,below=7mm of start.south,xshift=.235\textwidth] (geom)
 {\textbf{BSP}\\[3pt]
 Polynomial boundary tubes\\[3pt]
 Local-error peeling};
 \draw[flow] (start.south) -- ++(0,-3mm) -| (stat.north);
 \draw[flow] (start.south) -- ++(0,-3mm) -| (geom.north);
 \node[box,below=12mm of stat] (rough)
 {\textbf{Rough case: $\beta_-\le1$}\\[3pt]
 $\E B_r\lesssim r\Lambda_{\alpha,\beta_-}(\underline h)q_{\bm h}^{\alpha}$};
 \node[box,below=12mm of geom] (smooth)
 {\textbf{Smooth case: $\beta_->1$}\\[3pt]
 $\E B_r\lesssim r\Lambda_{\alpha,\beta_-}(\underline h)q_{\bm h}^{1+\alpha}
 +(r/\underline h)^2$};
 \coordinate (split) at ($(geom.south)+(0,-5mm)$);
 \draw[flow] (geom.south) -- (split) -| (rough.north);
 \draw[flow] (split) -- (smooth.north);
 \node[wide,below=9mm of smooth.south,xshift=-.235\textwidth] (end)
 {\textbf{Choose the bandwidths in \cref{cor:oracle-upper}}\\[3pt]
 $\E\cE_r\lesssim n^{-\frac{\bar\beta(1+\alpha)}{(2\bar\beta+d)}}
 +\psi_{\alpha,\beta_-}(r)$};
 \draw[flow] (rough.south) |- ($(end.north)+(0,4mm)$) -- (end.north);
 \draw[flow] (smooth.south) |- ($(end.north)+(0,4mm)$) -- (end.north);
 \draw[flow] (stat.west) -- ++(-5mm,0) |- (end.west);
\end{tikzpicture}
\caption{Proof sketch for the upper bound.}
\label{fig:upper-proof-roadmap}
\end{figure}
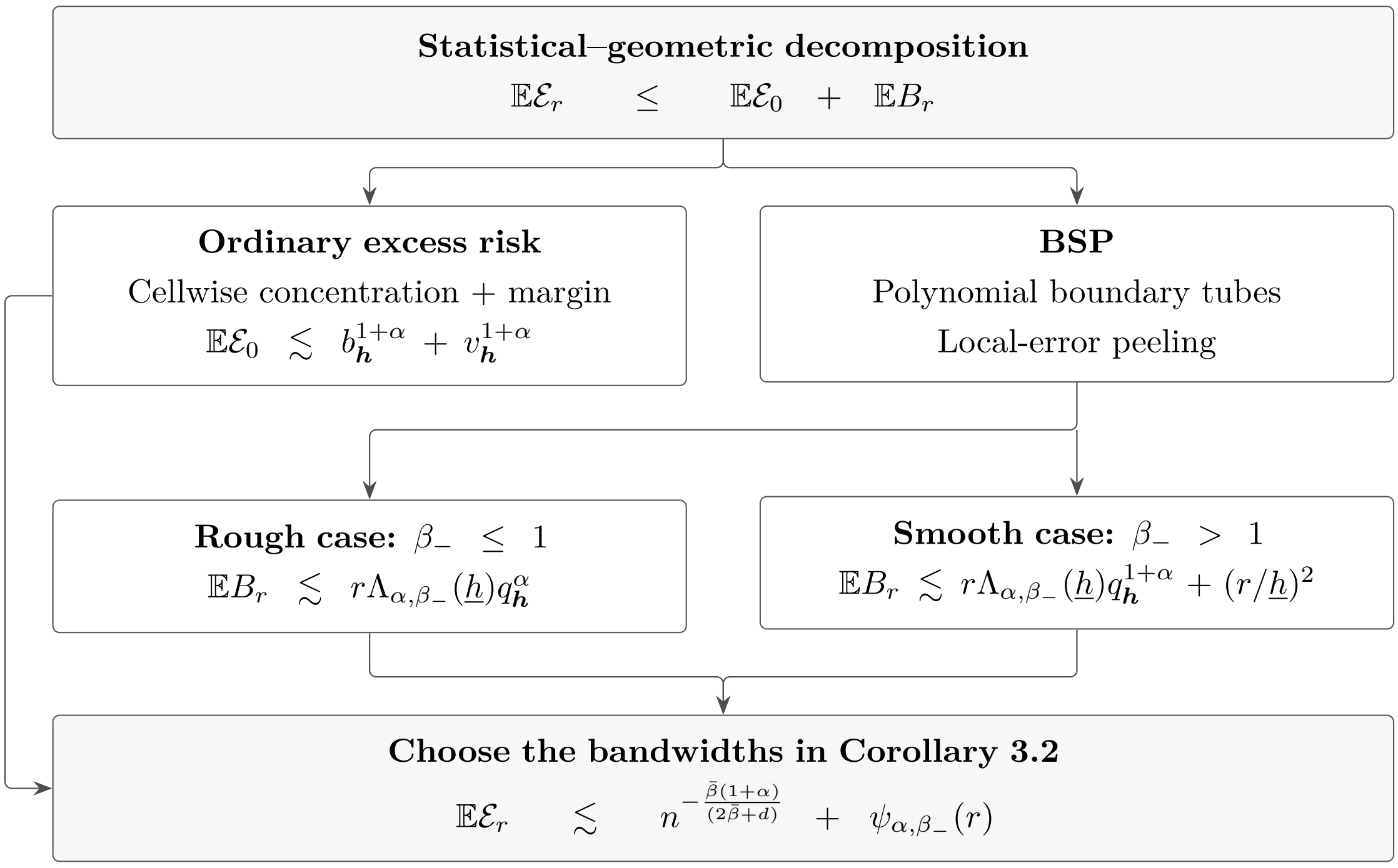

We first prove \cref{thm:fixed-upper}.
Let $v_{\bm h}:=(nV_{\bm h})^{-1/2}$ and
$q_{\bm h}:=1+(b_{\bm h}+v_{\bm h})/\underline h^{\beta_-}$.
By \cref{thm:decomposition},
\begin{equation}
\E\cE_r(\wh f_{\bm h})
\le\E\cE_0(\wh f_{\bm h})+\E B_r(\wh f_{\bm h}).
\label{eq:upper-decomposition}
\end{equation}
\Cref{fig:upper-proof-roadmap} summarizes the proof sketch.
The auxiliary lemmas are proved in \cref{app:upper-lemmas}.

\subsection{Ordinary excess risk}
\label{app:ordinary}

For $x\in Q$, $\wh f_{\bm h}(x)\ne f^*(x)$ implies
$\Delta(x)\le|\wh P_Q(x)-\eta(x)|$.
On $\{\Delta(x)>2C_{\rm app}b_{\bm h}\}$, applying
\cref{lem:cell-concentration} with $t=\Delta(x)/2$ gives
\[
P^{\otimes n}\{\wh f_{\bm h}(x)\ne f^*(x)\}
\le\1\{\Delta(x)\le2C_{\rm app}b_{\bm h}\}
 +Ce^{-c\Delta(x)^2/v_{\bm h}^2}.
\]
By the margin condition, we have
\begin{align*}
\E\cE_0(\wh f_{\bm h})
&=2\int_{[0,1]^d}\Delta(x)
 P^{\otimes n}\{\wh f_{\bm h}(x)\ne f^*(x)\}\,dP_X(x)\\
&\le Cb_{\bm h}^{1+\alpha}
 +C\E_X\{\Delta e^{-c\Delta^2/v_{\bm h}^2}\}.
\end{align*}
Splitting the last expectation over $\{\Delta\le v_{\bm h}\}$
and $\{2^kv_{\bm h}<\Delta\le2^{k+1}v_{\bm h}\}$, $k\ge0$,
and applying the margin condition yields
\begin{align*}
\E_X\{\Delta e^{-c\Delta^2/v_{\bm h}^2}\}
&\le v_{\bm h}P_X(\Delta\le v_{\bm h})
 +\sum_{k\ge0}2^{k+1}v_{\bm h}e^{-c4^k}
 P_X(\Delta\le2^{k+1}v_{\bm h})\\
&\le Cv_{\bm h}^{1+\alpha}
 \left\{1+\sum_{k\ge0}2^{(k+1)(1+\alpha)}e^{-c4^k}\right\}\\
&\le Cv_{\bm h}^{1+\alpha}.
\end{align*}
Consequently,
\begin{equation}
\E\cE_0(\wh f_{\bm h})
\le C\{b_{\bm h}^{1+\alpha}+v_{\bm h}^{1+\alpha}\}.
\label{eq:ordinary-final}
\end{equation}

\subsection{Boundary susceptibility probability}
\label{app:tube}

Assume $0<r\le\underline h$. Cover $[0,1]^d$ by the closed cubes
\[
S=\prod_{j=1}^d[k_j\underline h,(k_j+1)\underline h],
\qquad
k_j\in\{0,\ldots,\underline h^{-1}-1\}.
\]
These cubes have disjoint interiors and may cross the boundaries of
the rectangles $Q\in\cQ_{\bm h}$ used to construct $\wh f_{\bm h}$.
Hence each $S$ intersects $\overline Q$ for at most $3^d$ rectangles
$Q\in\cQ_{\bm h}$.
Set $\Gamma:=\partial_{[0,1]^d}\{\wh f_{\bm h}=1\}$.
For $A\subset\mathbb R^d$ and
$y=(y_k)_{k\ne j}\in\mathbb R^{d-1}$, let $N_j(A,y)$ denote the
number of points in the $j$th coordinate section of $A$ at $y$:
\[
N_j(A,y):=\#\left\{s\in\mathbb R:
(y_1,\ldots,y_{j-1},s,y_{j+1},\ldots,y_d)\in A\right\}.
\]
In the interior of each rectangle $Q\in\cQ_{\bm h}$, the boundary
is contained in $\{\wh P_Q=1/2\}$ and is empty if
$\wh P_Q\equiv1/2$. All other boundary points lie on rectangle faces.
Each cube $S$ meets at most $3^d$ rectangle closures, and the fitted
polynomials have degree at most $p$. Consequently, $\Gamma\cap S$
is compact and semialgebraic, with dimension at most $d-1$, and
\cref{lem:semialgebraic-tube} applies with a constant depending only
on $d$ and $p$.
The density upper bound then gives
\[
B_r(\wh f_{\bm h})
\le\mu_+\Leb(\Gamma^{+r})
\le\mu_+\sum_S\Leb\{(\Gamma\cap S)^{+r}\}.
\]
Applying \cref{lem:semialgebraic-tube} with $A=\Gamma\cap S$
and $t=\underline h$, then summing over $S$, gives
\begin{align}
B_r(\wh f_{\bm h})
&\le Cr\sum_{j=1}^d\int N_j(\Gamma,y)\,dy
 +C\1\{d\ge2\}(r/\underline h)^2,
\label{eq:tube-reduction}\\*
B_r(\wh f_{\bm h})
&\le Cr\underline h^{d-1}\#\{S:S\cap\Gamma\ne\varnothing\}.
\label{eq:active-tube}
\end{align}
We use \eqref{eq:tube-reduction} for $\beta_->1$ and
\eqref{eq:active-tube} for $\beta_-\le1$.

For \eqref{eq:tube-reduction}, each boundary point belongs to at most
$2^d$ cubes $S$, so
\[
\sum_S\int N_j(\Gamma\cap S,y)\,dy
\le2^d\int N_j(\Gamma,y)\,dy.
\]
Since there are $\underline h^{-d}$ cubes and $r\le\underline h$,
the local remainders satisfy
\[
\sum_S\sum_{k=2}^d r^k\underline h^{d-k}
=\sum_{k=2}^d(r/\underline h)^k
\le C(r/\underline h)^2.
\]
For \eqref{eq:active-tube}, almost every coordinate line meets
$\Gamma\cap S$ in at most $C_{d,p}$ points. The coordinate projection
of $S$ has $(d-1)$-dimensional measure $\underline h^{d-1}$, so
\[
\int N_j(\Gamma\cap S,y)\,dy\le C_{d,p}\underline h^{d-1}.
\]
Thus \cref{lem:semialgebraic-tube} and $r\le\underline h$ give
\[
\begin{aligned}
\Leb\{(\Gamma\cap S)^{+r}\}
&\le C\left(r\underline h^{d-1}
+\sum_{k=2}^d r^k\underline h^{d-k}\right)\\
&\le Cr\underline h^{d-1}.
\end{aligned}
\]
Summing over $S$ meeting $\Gamma$ proves \eqref{eq:active-tube}.

Define the local estimation errors and levels
\[
E_Q:=\min\{1,\|\wh P_Q-\eta\|_{\infty,Q}\},\quad
u_k:=2^ku_0,\quad u_0:=C_*(b_{\bm h}+v_{\bm h}),
\]
where $C_*>2C_{\rm app}$. Note that
\begin{equation}
z\in\Gamma\quad\Longrightarrow\quad
\Delta(z)\le\max_{Q:z\in\overline Q}E_Q.
\label{eq:boundary-error}
\end{equation}
Indeed, if $\Delta(z)>\max_{Q:z\in\overline Q}E_Q$, then each of
these errors is less than $1/2$.
By continuity, every adjacent fitted polynomial lies on
the same side of $1/2$ as $\eta(z)$ in a neighborhood of
$z$, contradicting $z\in\Gamma$.
For $u_k<1$, \cref{lem:cell-concentration} yields
\begin{equation}
P^{\otimes n}(E_Q>u_k)
\le Ce^{-cnV_{\bm h}(u_k-C_{\rm app}b_{\bm h})^2}
\le Ce^{-c4^k},
\label{eq:neighborhood-error}
\end{equation}
where $u_k-C_{\rm app}b_{\bm h}\ge u_k/2$ and
$nV_{\bm h}u_0^2\ge C_*^2$. 

\subsubsection{The regime \texorpdfstring{$\beta_-\le1$}{beta-min <= 1}}
\label{app:rough}

Define
\[
\mathcal S(u):=\{S:\min_{x\in S}\Delta(x)\le u\},\qquad
E_S:=\max_{Q:\overline Q\cap S\ne\varnothing}E_Q.
\]
If $S\in\mathcal S(u)$, let $y$ be a minimizer of $\Delta$ on $S$.
For every $x\in S$, we have $|x_j-y_j|\le\underline h$ for all $j$.
Thus, \cref{lem:anisotropic-consequences} gives
\begin{align*}
\Delta(x)
&\le\Delta(y)+|\eta(x)-\eta(y)|\\
&\le u+L\sum_{j=1}^d|x_j-y_j|^{1\wedge\beta_j}\\
&\le u+dL\underline h^{\beta_-},
\end{align*}
where the last inequality uses $\beta_-\le1\wedge\beta_j$ and
$0<\underline h\le1$. Consequently,
\[
S\in\mathcal S(u)\quad\Longrightarrow\quad
S\subseteq\{\Delta\le u+dL\underline h^{\beta_-}\}.
\]
By the density lower bound and the margin condition, we have
\[
\underline h^d\#\mathcal S(u)
=\Leb\!\left(\bigcup_{S\in\mathcal S(u)}S\right)
\le\mu_-^{-1}P_X(\Delta\le u+C\underline h^{\beta_-})
\le C(u+\underline h^{\beta_-})^\alpha.
\]
Dividing by $\underline h$ gives
\begin{equation}
\underline h^{d-1}\#\mathcal S(u)
\le C\underline h^{-1}(u+\underline h^{\beta_-})^\alpha
=C\Lambda_{\alpha,\beta_-}(\underline h)
 (1+u/\underline h^{\beta_-})^\alpha.
\label{eq:rough-count}
\end{equation}
If $S\cap\Gamma\ne\varnothing$, choose $z\in S\cap\Gamma$.
Thus, by \eqref{eq:boundary-error},
\[
\min_{x\in S}\Delta(x)
\le\Delta(z)
\le\max_{Q:z\in\overline Q}E_Q
\le E_S,
\]
which implies $S\in\mathcal S(E_S)$. Splitting according to $E_S\le u_0$ and
$u_k<E_S\le u_{k+1}$ gives
\begin{align*}
\1\{S\cap\Gamma\ne\varnothing\}
&=\1\{S\cap\Gamma\ne\varnothing\}\1\{E_S\le u_0\}+\sum_{k\ge0}\1\{S\cap\Gamma\ne\varnothing\}
 \1\{u_k<E_S\le u_{k+1}\}\\*
&\le\1\{S\in\mathcal S(u_0)\}
 +\sum_{k\ge0}\1\{S\in\mathcal S(u_{k+1})\}\1\{E_S>u_k\}.
\end{align*}
At most $3^d$ cells $Q\in\mathcal Q_{\bm h}$ satisfy
$\overline Q\cap S\ne\varnothing$, so a union bound gives
\[
P^{\otimes n}(E_S>u_k)
\le\sum_{Q:\overline Q\cap S\ne\varnothing}
 P^{\otimes n}(E_Q>u_k)
\le Ce^{-c4^k}.
\]
Taking expectations and summing over $S$, we obtain
\begin{align}
\E\#\{S:S\cap\Gamma\ne\varnothing\}
&\le\#\mathcal S(u_0)
 +\sum_{k\ge0}\sum_{S\in\mathcal S(u_{k+1})}
 P^{\otimes n}(E_S>u_k)\notag\\
&\le\#\mathcal S(u_0)
 +C\sum_{k\ge0}e^{-c4^k}\#\mathcal S(u_{k+1}).
\label{eq:peeling}
\end{align}
Since $1+u_{k+1}/\underline h^{\beta_-}\le C2^kq_{\bm h}$,
\eqref{eq:rough-count} implies
\[
\underline h^{d-1}\#\mathcal S(u_{k+1})
\le C\Lambda_{\alpha,\beta_-}(\underline h)
 2^{k\alpha}q_{\bm h}^{\alpha}.
\]
Combining this with \eqref{eq:active-tube} and \eqref{eq:peeling} yields
\begin{equation}
\E B_r(\wh f_{\bm h})
\le Cr\Lambda_{\alpha,\beta_-}(\underline h)q_{\bm h}^{\alpha}
\left(1+\sum_{k\ge0}e^{-c4^k}2^{k\alpha}\right)
\le Cr\Lambda_{\alpha,\beta_-}(\underline h)q_{\bm h}^{\alpha}.
\label{eq:rough-G-sharp}
\end{equation}

\subsubsection{The regime \texorpdfstring{$\beta_->1$}{beta-min > 1}}
\label{app:smooth}
To control the boundary intersections in \eqref{eq:tube-reduction},
we first bound the number of cell sections where $\eta$ is close
to $1/2$. We then combine this bound with the local estimation
error bounds.

For $y=(y_k)_{k\ne j}\in[0,1]^{d-1}$, define
\[
g_y(s):=\eta(y_1,\ldots,y_{j-1},s,y_{j+1},\ldots,y_d)-\tfrac12,
\qquad s\in[0,1].
\]
For each $Q\in\mathcal Q_{\bm h}$, let $Q_y$ denote the closure of
its $j$th coordinate section at $y$:
\[
Q_y:=\overline{\{s\in[0,1]:
(y_1,\ldots,y_{j-1},s,y_{j+1},\ldots,y_d)\in Q\}}.
\]
Let $K_j(y,u)$ count the sections meeting $\{|g_y|\le u\}$:
\[
K_j(y,u):=\#\bigl\{Q\in\mathcal Q_{\bm h}:
Q_y\cap\{|g_y|\le u\}\ne\varnothing\bigr\}.
\]
Every section counted by $K_j(y,u)$ contains a point where $|g_y|\le u$.
Since its length is $h_j$,
\[
\bigcup_{\substack{Q\in\mathcal Q_{\bm h}\\
Q_y\cap\{|g_y|\le u\}\ne\varnothing}} Q_y
\subseteq[0,1]\cap\{|g_y|\le u\}^{+h_j}.
\]
Let $T=h_j+u^{1/\beta_-}$. If $T\le1$, then $h_j\le T$ and
$u\le T^{\beta_-}$, so the last set is contained in
$[0,1]\cap\{|g_y|\le T^{\beta_-}\}^{+T}$.
\cref{lem:axis-tube} gives
\begin{align*}
h_j\int K_j(y,u)\,dy
&\le\int\Leb_1\bigl([0,1]\cap\{|g_y|\le u\}^{+h_j}\bigr)\,dy\\
&\le CT\Lambda_{\alpha,\beta_-}(T).
\end{align*}
Also, $T\ge h_j$ gives
$\Lambda_{\alpha,\beta_-}(T)\le\Lambda_{\alpha,\beta_-}(h_j)$, and
$T/h_j\le2(1+u/h_j^{\beta_-})^{1/\beta_-}$.
Dividing by $h_j$ therefore gives the first inequality in
\eqref{eq:axis-count} when $T\le1$.

Each coordinate line meets at most $h_j^{-1}$ cells, independently of $T$.
Thus $K_j(y,u)\le h_j^{-1}$ for every $y$ and $u\ge0$.
For $T>1$, we use this bound to obtain
$\int K_j(y,u)\,dy\le h_j^{-1}$.
Moreover, $T>1$ implies $\max\{h_j,u^{1/\beta_-}\}\ge1/2$, and hence
\[
\left(1+\frac{u}{h_j^{\beta_-}}\right)^{1/\beta_-}
\ge\frac{\max\{h_j,u^{1/\beta_-}\}}{h_j}
\ge\frac1{2h_j}.
\]
Since $\Lambda_{\alpha,\beta_-}(h_j)\ge1$, this also gives the first
inequality in \eqref{eq:axis-count}, after enlarging $C$.
Combining the two cases and using $\underline h\le h_j$ and
$1/\beta_-<1+\alpha$, we obtain, for all $u\ge0$,
\begin{equation}
\label{eq:axis-count}
\begin{split}
\int K_j(y,u)\,dy
&\le C\Lambda_{\alpha,\beta_-}(h_j)
 (1+u/h_j^{\beta_-})^{1/\beta_-}\\
&\le C\Lambda_{\alpha,\beta_-}(\underline h)
 (1+u/\underline h^{\beta_-})^{1+\alpha}.
\end{split}
\end{equation}
On almost every coordinate line, each nonzero threshold polynomial
has at most $p$ roots per cell. Assign each boundary point to an
incident cell maximizing $E_Q$. Each cell receives at most $p$ interior
points and two endpoints, and every receiving cell satisfies
$\min_{s\in Q_y}|g_y(s)|\le E_Q$ by \eqref{eq:boundary-error}.
Splitting according to $E_Q\le u_0$ or $u_k<E_Q\le u_{k+1}$ gives
\[
\begin{aligned}
\frac{N_j(\Gamma,y)}{p+2}
&\le\sum_{\substack{Q\in\mathcal Q_{\bm h}\\Q_y\ne\varnothing}}
 \1\left\{\min_{s\in Q_y}|g_y(s)|\le E_Q\right\}\\
&\le K_j(y,u_0)
+\sum_{k\ge0}\sum_{Q\in\mathcal Q_{\bm h}}
 \1\{Q_y\cap\{|g_y|\le u_{k+1}\}\ne\varnothing\}\1\{E_Q>u_k\}.
\end{aligned}
\]
For fixed $y$, the indicators involving $g_y$ are deterministic,
so \eqref{eq:neighborhood-error} gives
\[
\begin{split}
&\E\sum_{Q\in\mathcal Q_{\bm h}}
 \1\{Q_y\cap\{|g_y|\le u_{k+1}\}\ne\varnothing\}\1\{E_Q>u_k\}\\
=&\sum_{Q\in\mathcal Q_{\bm h}}
 \1\{Q_y\cap\{|g_y|\le u_{k+1}\}\ne\varnothing\}
 P^{\otimes n}(E_Q>u_k)\\
\le & Ce^{-c4^k}K_j(y,u_{k+1}).
\end{split}
\]
Note that $1+u_{k+1}/\underline h^{\beta_-}\le C2^kq_{\bm h}$. Tonelli's theorem and \eqref{eq:axis-count} give
\begin{align*}
\E\int N_j(\Gamma,y)\,dy
&\le C\int K_j(y,u_0)\,dy
+C\sum_{k\ge0}e^{-c4^k}\int K_j(y,u_{k+1})\,dy\\*
&\le C\Lambda_{\alpha,\beta_-}(\underline h)q_{\bm h}^{1+\alpha}
\left(1+\sum_{k\ge0}e^{-c4^k}2^{k(1+\alpha)}\right).
\end{align*}
Substitution into \eqref{eq:tube-reduction} gives
\begin{equation}
\E B_r(\wh f_{\bm h})
\le Cr\Lambda_{\alpha,\beta_-}(\underline h)q_{\bm h}^{1+\alpha}
+C\1\{d\ge2\}(r/\underline h)^2.
\label{eq:smooth-G}
\end{equation}

\subsection{Combination of the upper bound}
\label{app:upper-completion}
Since $q_{\bm h}\ge1$, the factor $q_{\bm h}^\alpha$ in
\eqref{eq:rough-G-sharp} is at most $q_{\bm h}^{1+\alpha}$.
Thus \eqref{eq:upper-decomposition}, \eqref{eq:ordinary-final},
\eqref{eq:rough-G-sharp}, and \eqref{eq:smooth-G}, together with
$B_0=0$, prove \cref{thm:fixed-upper}.

\Needspace{10\baselineskip}
For \cref{cor:oracle-upper}, choose each $h_j$ as the reciprocal
of a positive integer so that
\[
\begin{split}
n^{-\frac{\bar\beta}{\beta_j(2\bar\beta+d)}}
\vee r^{\frac{\beta_-}{\beta_j(\beta_-+1)}}
\le h_j\le 2\left(
n^{-\frac{\bar\beta}{\beta_j(2\bar\beta+d)}}
\vee r^{\frac{\beta_-}{\beta_j(\beta_-+1)}}
\right),\qquad j=1,\ldots,d.
\end{split}
\]
These bandwidths satisfy
\[
nV_{\bm h}\ge n^{2\bar\beta/(2\bar\beta+d)}\ge1,
\qquad
\underline h\ge r^{1/(\beta_-+1)}\ge r.
\]
Thus \cref{thm:fixed-upper} applies. Moreover,
\begin{align*}
v_{\bm h}&\le n^{-\bar\beta/(2\bar\beta+d)},\\
b_{\bm h}+v_{\bm h}
&\le C\left\{
n^{-\bar\beta/(2\bar\beta+d)}
\vee r^{\beta_-/(\beta_-+1)}\right\}
\le C\underline h^{\beta_-}.
\end{align*}
Hence $q_{\bm h}\le C$, and \eqref{eq:ordinary-final} gives
\begin{equation}
\E\cE_0(\wh f_{\bm h})
\le C\left\{
n^{-\bar\beta(1+\alpha)/(2\bar\beta+d)}
+r^{\beta_-(1+\alpha)/(\beta_-+1)}\right\}.
\label{eq:ordinary-selected-bandwidth}
\end{equation}

By \eqref{eq:Lambda}, \eqref{eq:psi-definition}, and
$\underline h\ge r^{1/(\beta_-+1)}$,
\begin{align*}
r\Lambda_{\alpha,\beta_-}(\underline h)
&\le r\Lambda_{\alpha,\beta_-}\bigl(r^{1/(\beta_-+1)}\bigr)
\le C\psi_{\alpha,\beta_-}(r),\\
r^{\beta_-(1+\alpha)/(\beta_-+1)}+r
&\le C\psi_{\alpha,\beta_-}(r).
\end{align*}
The remainder in \eqref{eq:smooth-G} occurs only when $\beta_->1$,
in which case
\[
(r/\underline h)^2\le r^{2\beta_-/(\beta_-+1)}\le r.
\]
Thus \eqref{eq:rough-G-sharp} and \eqref{eq:smooth-G} imply
$\E B_r(\wh f_{\bm h})\le C\psi_{\alpha,\beta_-}(r)$.
Combining this with \eqref{eq:upper-decomposition} and
\eqref{eq:ordinary-selected-bandwidth} yields
\[
\sup_{P\in\cP(\alpha,\bm\beta,L)}
\E\cE_r(\wh f_{\bm h})
\le C\left\{n^{-\bar\beta(1+\alpha)/(2\bar\beta+d)}
+\psi_{\alpha,\beta_-}(r)\right\}.
\]
This proves \cref{cor:oracle-upper}.

\subsection{Auxiliary technical lemmas}
\label{app:upper-lemmas}

For a multi-index $\bm k\in\mathbb N_0^d$, write
\[
D^{\bm k}:=
\frac{\partial^{|\bm k|}}
{\partial x_1^{k_1}\cdots\partial x_d^{k_d}}.
\]

\begin{lemma}
\label{lem:anisotropic-consequences}
Under \cref{ass:holder}, for all $x,y\in[0,1]^d$,
\begin{equation}
|\eta(x)-\eta(y)|
\le
L\sum_{j=1}^d|x_j-y_j|^{1\wedge\beta_j}
\le
dL\|x-y\|_2^{1\wedge\beta_-}.
\label{eq:coordinate-modulus}
\end{equation}

Moreover, \cref{ass:holder} holds with each $\beta_j$ replaced by $\beta_-$, that is,
\begin{equation}
\eta\in\mathcal C^{(\beta_-,\ldots,\beta_-)}(L).
\label{eq:coordinate-section-regularity}
\end{equation}
Fix an integer $p\ge1$ such that $p+1>\beta_+$.
For each rectangle $Q\subset[0,1]^d$ with lengths $\bm h=(h_1,\ldots,h_d)$, there exists a polynomial
$P_Q$ of total degree at most $p$ such that
\begin{equation}
\|\eta-P_Q\|_{\infty,Q}
\le
C\sum_{j=1}^d h_j^{\beta_j}.
\label{eq:coordinate-polynomial-bias}
\end{equation}
Here $C$ depends only on $d,\bm\beta,L$ and $p$.
\end{lemma}

\begin{proof}
The H\"older bound for $\beta_j\le1$ and the mean value theorem for
$\beta_j>1$ give
\[
|\eta(x+te_j)-\eta(x)|\le L|t|^{1\wedge\beta_j},
\qquad x,x+te_j\in[0,1]^d.
\]
Changing one coordinate at a time yields
\begin{align*}
|\eta(x)-\eta(y)|
&\le\sum_{j=1}^d\bigl|
\eta(x_1,\ldots,x_j,y_{j+1},\ldots,y_d)
-\eta(x_1,\ldots,x_{j-1},y_j,\ldots,y_d)\bigr|\\
&\le L\sum_{j=1}^d|x_j-y_j|^{1\wedge\beta_j}
\le dL\|x-y\|_2^{1\wedge\beta_-},
\end{align*}
which proves \eqref{eq:coordinate-modulus}.
For \eqref{eq:coordinate-section-regularity}, let
$m:=\lceil\beta_-\rceil-1$ and recall $m_j=\lceil\beta_j\rceil-1$.
If $m_j=m$, lower the H\"older exponent
from $\beta_j-m$ to $\beta_--m$. If $m_j>m$, the mean value theorem
gives a Lipschitz bound with constant $L$ for $\partial_j^m\eta$,
which also gives the required H\"older bound on $[0,1]$.

For \eqref{eq:coordinate-polynomial-bias}, let
\[
N=\left\lceil\log_2\frac1{\max_jh_j}\right\rceil,
\qquad R_\ell=(z_Q+2^\ell(Q-z_Q))\cap[0,1]^d,
\quad 0\le\ell\le N.
\]
The $j$th side length of $R_\ell$ lies in
$[2^{\ell-1}h_j,2^\ell h_j]$.
Let $U_\ell$ be the tensor Lagrange interpolant of $\eta$ on $R_\ell$, with
degree at most $m_j$ in coordinate $j$ and fixed rescaled nodes.
Stability of interpolation and the univariate Taylor bounds give
\[
\|U_\ell\|_{\infty,R_\ell}\le C,
\qquad
\|\eta-U_\ell\|_{\infty,R_\ell}
\le C\sum_j2^{\ell\beta_j}h_j^{\beta_j}.
\]
For $|\bm k|=p+1$ and
$1\le\ell\le N$, rescaling $R_{\ell-1}$ to $[0,1]^d$ and applying
the polynomial derivative bound gives
\begin{align*}
\left(\prod_jh_j^{k_j}\right)
\|D^{\bm k}(U_{\ell-1}-U_\ell)\|_{\infty,Q}
&\le C2^{-\ell(p+1)}
\|U_{\ell-1}-U_\ell\|_{\infty,R_{\ell-1}}\\
&\le C\sum_j2^{-\ell(p+1-\beta_j)}h_j^{\beta_j}.
\end{align*}
Since $U_0=U_N+\sum_{\ell=1}^N(U_{\ell-1}-U_\ell)$, the same
derivative bound gives
\begin{align*}
\left(\prod_jh_j^{k_j}\right)
\|D^{\bm k}U_0\|_{\infty,Q}
&\le C2^{-N(p+1)}
+C\sum_jh_j^{\beta_j}\sum_{\ell=1}^N2^{-\ell(p+1-\beta_j)}\\
&\le Cb_{\bm h},
\end{align*}
since $2^{-N(p+1)}\le(\max_jh_j)^{p+1}\le b_{\bm h}$ and
$p+1>\beta_+$.
Let $P_Q$ be the Taylor polynomial of order $p$ of $U_0$ at $z_Q$.
Taylor's theorem and the preceding derivative estimate give
\begin{align*}
\|U_0-P_Q\|_{\infty,Q}
&\le C\sum_{|\bm k|=p+1}\left(\prod_jh_j^{k_j}\right)
\|D^{\bm k}U_0\|_{\infty,Q}\\
&\le Cb_{\bm h}.
\end{align*}
Combining this with the interpolation bound yields
\[
\|\eta-P_Q\|_{\infty,Q}
\le\|\eta-U_0\|_{\infty,Q}+\|U_0-P_Q\|_{\infty,Q}
\le Cb_{\bm h}.
\]
\end{proof}

\begin{lemma}
\label{lem:cell-concentration}
There exist $C_{\rm app},c,C>0$ such that, for
$nV_{\bm h}\ge1$, $Q\in\cQ_{\bm h}$, and $0<t\le1$,
\begin{equation}
\sup_{P\in\cP(\alpha,\bm\beta,L)}
P^{\otimes n}\{\|\wh P_Q-\eta\|_{\infty,Q}>C_{\rm app}b_{\bm h}+t\}
\le Ce^{-cnV_{\bm h}t^2}.
\label{eq:cell-concentration}
\end{equation}
\end{lemma}
\begin{proof}
Let $P_Q$ be the polynomial in \cref{lem:anisotropic-consequences}.
Bernstein's inequality for the Gram entries gives
\[
\mathcal A_Q:=\{\|\wh G_Q-\E\wh G_Q\|_{\rm op}\le\kappa_p/2\},
\qquad P^{\otimes n}(\mathcal A_Q^c)\le Ce^{-cnV_{\bm h}}.
\]
On $\mathcal A_Q$, 
\[
\lambda_{\min}(\wh G_Q)\ge\kappa_p/2,\qquad
\|\wh G_Q^{-1}\|_{\rm op}\le2/\kappa_p.
\]
The constant basis function gives
\[
\frac1{nV_{\bm h}}\sum_i\1\{X_i\in Q\}
=(\wh G_Q)_{00}\le\|\wh G_Q\|_{\rm op}\le C.
\]
On the same event, for $x\in Q$,
\[
\wh P_Q(x)-P_Q(x)
=\phi_Q(x)^\top\wh G_Q^{-1}(\xi_Q+R_Q),
\]
where
\begin{align*}
\xi_Q&:=\frac1{nV_{\bm h}}\sum_i\1\{X_i\in Q\}
 \phi_Q(X_i)\{Y_i-\eta(X_i)\},\\
R_Q&:=\frac1{nV_{\bm h}}\sum_i\1\{X_i\in Q\}
 \phi_Q(X_i)\{\eta(X_i)-P_Q(X_i)\}.
\end{align*}
By \eqref{eq:coordinate-polynomial-bias} and the preceding count,
$\|R_Q\|_2\le Cb_{\bm h}$ on $\mathcal A_Q$.
Bernstein's inequality also gives
\[
P^{\otimes n}(\|\xi_Q\|_2>s)\le Ce^{-cnV_{\bm h}s^2},
\qquad 0<s\le1.
\]
Consequently, on $\mathcal A_Q$,
\[
\|\wh P_Q-\eta\|_{\infty,Q}
\le C(\|\xi_Q\|_2+b_{\bm h}).
\]
Taking $s=ct$ gives
\begin{align*}
&P^{\otimes n}\{\|\wh P_Q-\eta\|_{\infty,Q}
 >C_{\rm app}b_{\bm h}+t\}\\
&\quad\le P^{\otimes n}(\mathcal A_Q^c)
 +P^{\otimes n}(\|\xi_Q\|_2>ct)\\
&\quad\le Ce^{-cnV_{\bm h}t^2}.
\end{align*}
\end{proof}

\Needspace{10\baselineskip}
\begin{lemma}\label{lem:semialgebraic-tube}
Let $A\subset\R^d$ be compact, semialgebraic, of dimension at most
$d-1$, and contained in a cube of side $t$. Suppose $A$ is defined by a bounded number of polynomial
equalities and inequalities of bounded degree. Then
\begin{equation}
\Leb(A^{+r})\le Cr\sum_{j=1}^d\int N_j(A,y)\,dy
+C\sum_{k=2}^d r^kt^{d-k},\qquad 0<r\le t.
\label{eq:semialgebraic-tube}
\end{equation}
The constant $C$ depends only on
$d$ and the bounds on the number and degrees of these polynomials.
\end{lemma}
\begin{proof}
For $d=1$, $\Leb(A^{+r})\le2r\#A$.
For $d\ge2$, take the grid with vertices $u+rz$,
$z\in\mathbb Z^d$, where $u$ is uniform on $[0,r)^d$.
Let $H_r$ count its closed cubes of side $r$ meeting $A$.
These cubes cover $A$. Extending each by $r$ on every side gives
$H_r$ cubes of side $3r$ covering $A^{+r}$. Hence
\[
\Leb(A^{+r})\le(3r)^dH_r.
\]
Since $\Leb(A)=0$, grid vertices miss $A$ almost surely.
If $A$ misses a counted cube's boundary, the cube contains
a whole component of $A$.
Otherwise, repeat the argument on a proper face meeting $A$,
within its coordinate plane.
Each face belongs to at most $2^d$ cubes.

For each $j$, let $\mathcal L_j$ be the family of grid lines
\[
\Pi=\{x\in\R^d:x_\ell=u_\ell+rz_\ell,\ \ell\in\{1,\ldots,d\}\setminus\{j\}\},
\qquad (z_\ell)_{\ell\ne j}\in\mathbb Z^{d-1}.
\]
For $2\le k\le d$, each $k$-dimensional coordinate section has at most $C$ components
\citep{FriedlandYomdin2014}.
Fixing $d-k$ grid coordinates gives at most
$C(1+t/r)^{d-k}\le C(t/r)^{d-k}$ such planes meeting the cube containing $A$.
Thus
\[
H_r\le C\sum_{j=1}^d\sum_{\Pi\in\mathcal L_j}\#(A\cap\Pi)
+C\sum_{k=2}^d(t/r)^{d-k}.
\]
Averaging over the grid shift gives
\[
\begin{aligned}
\E\sum_{\Pi\in\mathcal L_j}\#(A\cap\Pi)
&=\frac1{r^{d-1}}\int_{[0,r)^{d-1}}
  \sum_{z\in\mathbb Z^{d-1}}N_j(A,y+rz)\,dy\\*
&=\frac1{r^{d-1}}\int_{\mathbb R^{d-1}}N_j(A,y)\,dy.
\end{aligned}
\]
Since $A$ is fixed, averaging the covering bound over $u$ gives
\[
\begin{aligned}
\Leb(A^{+r})
&\le3^dr^d\E H_r\\*
&\le Cr\sum_{j=1}^d\int N_j(A,y)\,dy
  +Cr^d\sum_{k=2}^d(t/r)^{d-k}\\*
&=Cr\sum_{j=1}^d\int N_j(A,y)\,dy
  +C\sum_{k=2}^dr^kt^{d-k}.
\end{aligned}
\]
\end{proof}

\begin{lemma}
\label{lem:axis-tube}
Suppose $\beta_->1$. Fix $j\in\{1,\ldots,d\}$ and let $g_y$
be the coordinate section defined in \cref{app:smooth}.
Define
$A_t(y)=\{s\in[0,1]:|g_y(s)|\le t^{\beta_-}\}$.
Under \cref{ass:holder,ass:margin,ass:density},
\begin{equation}
\int\Leb_1([0,1]\cap A_t(y)^{+t})\,dy
\le Ct\Lambda_{\alpha,\beta_-}(t),\qquad 0<t\le1.
\label{eq:axis-tube}
\end{equation}
\end{lemma}
\begin{proof}
Let $m:=\lceil\beta_-\rceil-1$.
By \eqref{eq:coordinate-section-regularity}, uniformly in $y$,
\[
\begin{gathered}
\max_{0\le a\le m}\|g_y^{(a)}\|_\infty\le C,\\
|g_y^{(m)}(s)-g_y^{(m)}(t)|
\le L|s-t|^{\beta_--m},\qquad s,t\in[0,1].
\end{gathered}
\]
Fix $y$ and let $L_t=\{z:|g_y(z)|=t^{\beta_-}\}$. Continuity gives
\[
[0,1]\cap A_t(y)^{+t}\subseteq A_t(y)\cup([0,1]\cap L_t^{+t}).
\]
Choose a fixed constant $A$ large enough that
\[
\rho(z):=\max_{1\le a\le m}
 (|g_y^{(a)}(z)|/A)^{1/(\beta_--a)}\le1/2.
\]

If $z\in L_t$ and $\rho(z)\le t$, Taylor expansion gives,
for $|s-z|\le t$,
\[
|g_y(s)|\le t^{\beta_-}
+A\sum_{a=1}^m\frac{t^{\beta_--a}|s-z|^a}{a!}
+C|s-z|^{\beta_-}\le Ct^{\beta_-}.
\]
The $t$-neighborhood of these level points lies in
$\{|g_y|\le Ct^{\beta_-}\}$.

For $\rho(z)>t$, define layers at scales
$\ell=2^kt$, $k=1,\ldots,\lfloor\log_2(1/t)\rfloor$, by
\[
L_\ell:=\{z\in L_t:\ell/2<\rho(z)\le\ell\}.
\]
At each $z\in L_\ell$,
\[
|g_y^{(b)}(z)|\le A\ell^{\beta_--b}\quad(1\le b\le m),\qquad
|g_y^{(a)}(z)|>A(\ell/2)^{\beta_--a}\quad\text{for some }a.
\]
Choose $c>0$ small enough for all $a=1,\ldots,m$.
For $s\in[z-c\ell,z+c\ell]\cap[0,1]$, Taylor expansion gives
\begin{align*}
|g_y^{(a)}(s)-g_y^{(a)}(z)|
&\le\sum_{b=a+1}^m
 \frac{|g_y^{(b)}(z)|}{(b-a)!}|s-z|^{b-a}
 +C|s-z|^{\beta_--a}\\
&\le\left(A\sum_{b=a+1}^m c^{b-a}+Cc^{\beta_--a}\right)
 \ell^{\beta_--a}\\
&<\tfrac12A(\ell/2)^{\beta_--a}.
\end{align*}
Thus $g_y^{(a)}$ has constant nonzero sign on this interval.
Rolle's theorem gives at most $a$ roots of each equation
$g_y=\pm t^{\beta_-}$ on the interval.

Choose a maximal $c\ell/8$-separated set of centers
$z_1,\ldots,z_{M_\ell}$ in $L_\ell$.
Taylor expansion at each center gives
\[
\sup_{|s-z_i|\le c\ell/16}|g_y(s)|
\le t^{\beta_-}+A\sum_{b=1}^m\ell^{\beta_--b}
 \frac{(c\ell/16)^b}{b!}+C(c\ell)^{\beta_-}
\le C\ell^{\beta_-}.
\]
Each truncated interval has length at least $c\ell/16$, so
\[
M_\ell\ell\le C\Leb_1\{|g_y|\le C\ell^{\beta_-}\}.
\]
By maximality, the intervals of radius $c\ell/8$ cover $L_\ell$.
The root bound gives at most $2m$ points per interval, so
\[
\#L_\ell\le2mM_\ell.
\]
Covering these points by intervals of radius $t$ gives
\[
\Leb_1([0,1]\cap L_\ell^{+t})
\le2t\,\#L_\ell
\le4mtM_\ell
\le Ct\ell^{-1}\Leb_1\{|g_y|\le C\ell^{\beta_-}\}.
\]
Combining this with the bound for $\rho(z)\le t$ gives
\begin{equation}
\Leb_1([0,1]\cap A_t(y)^{+t})
\le\Leb_1\{|g_y|\le Ct^{\beta_-}\}
+Ct\sum_\ell\ell^{-1}\Leb_1\{|g_y|\le C\ell^{\beta_-}\}.
\label{eq:fiber-tube-bound}
\end{equation}
By the density and margin assumptions,
\[
\int\Leb_1\{|g_y|\le C\ell^{\beta_-}\}\,dy
\le\mu_-^{-1}P_X(\Delta\le C\ell^{\beta_-})
\le C\ell^{\alpha\beta_-}.
\]
Summing the geometric series and using \eqref{eq:Lambda}, we obtain
\[
\int\Leb_1([0,1]\cap A_t(y)^{+t})\,dy
\le Ct^{\alpha\beta_-}
+Ct\sum_{k=1}^{\lfloor\log_2(1/t)\rfloor}
 (2^kt)^{\alpha\beta_--1}
\le Ct\Lambda_{\alpha,\beta_-}(t).
\]
\end{proof}

\Needspace{14\baselineskip}

\section{Proof of the lower bound}
\label{app:lower}
By \cref{thm:decomposition}, every classifier $f$ satisfies
\begin{equation}
\cE_r(f)
\ge
\frac12\{\cE_0(f)+B_r(f)\}.
\label{eq:lower-decomposition}
\end{equation}
We obtain the statistical and adversarial lower bounds
from separate constructions, each with $P_X$ uniform on
$[0,1]^d$.

\subsection{Statistical component}
\label{app:lower-statistical}

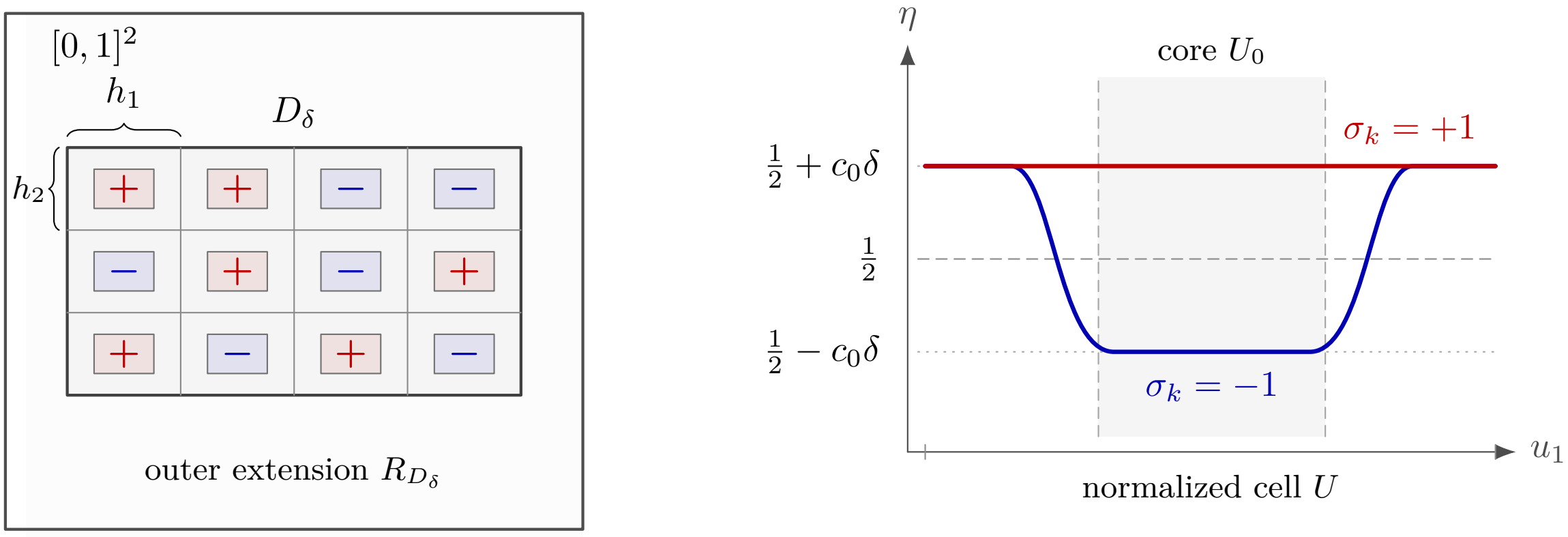
\begin{figure}[H]
\centering
\begin{tikzpicture}[
  x=1cm,
  y=1cm,
  font=\small,
  line cap=round,
  line join=round,
  >={Latex[length=2mm]}
]
  \fill[black!1] (0,0) rectangle (5.4,5.1);
  \draw[thick,black!70] (-0.2,0.1) rectangle (5.4,5.1);
  \node[anchor=north west] at (0.10,5.10) {$[0,1]^2$};

  \fill[black!4] (0.40,1.40) rectangle (4.80,3.80);
  \draw[thick,black!75] (0.40,1.40) rectangle (4.80,3.80);
  \node[anchor=south] at (2.60,3.86) {$D_\delta$};
  \foreach \xgrid in {1.50,2.60,3.70}
    \draw[black!45] (\xgrid,1.40) -- (\xgrid,3.80);
  \foreach \ygrid in {2.20,3.00}
    \draw[black!45] (0.40,\ygrid) -- (4.80,\ygrid);

  \draw[decorate,decoration={brace,amplitude=3.5pt}]
    (0.40,3.91) -- (1.50,3.91)
    node[midway,above=4pt] {$h_1$};
  \draw[decorate,decoration={brace,amplitude=2.8pt}]
    (0.32,3.00) -- (0.32,3.80);
  \node[font=\footnotesize] at (0.03,3.40) {$h_2$};

  \foreach \cx/\cy in {
    0.95/1.80,3.15/1.80,
    2.05/2.60,4.25/2.60,
    0.95/3.40,2.05/3.40
  }{
    \fill[red!75!black,opacity=.08]
      ($(\cx,\cy)+(-.29,-.19)$) rectangle ($(\cx,\cy)+(.29,.19)$);
    \draw[black!55]
      ($(\cx,\cy)+(-.29,-.19)$) rectangle ($(\cx,\cy)+(.29,.19)$);
    \node[text=red!75!black] at (\cx,\cy) {$\boldsymbol{+}$};
  }
  \foreach \cx/\cy in {
    2.05/1.80,4.25/1.80,
    0.95/2.60,3.15/2.60,
    3.15/3.40,4.25/3.40
  }{
    \fill[blue!70!black,opacity=.08]
      ($(\cx,\cy)+(-.29,-.19)$) rectangle ($(\cx,\cy)+(.29,.19)$);
    \draw[black!55]
      ($(\cx,\cy)+(-.29,-.19)$) rectangle ($(\cx,\cy)+(.29,.19)$);
    \node[text=blue!70!black] at (\cx,\cy) {$\boldsymbol{-}$};
  }
  \node[font=\footnotesize] at (2.60,0.65)
    {outer extension $R_{D_\delta}$};
  \node[font=\small] at (2.60,-0.38)
    {(a) Anisotropic tiling};

  \draw[->,black!70] (8.55,0.85) -- (8.55,4.80)
    node[above] {$\eta$};
  \draw[->,black!70] (8.55,0.85) -- (14.45,0.85)
    node[right] {$u_1$};

  \fill[black!4] (10.40,1.00) rectangle (12.60,4.48);
  \draw[densely dashed,black!35] (10.40,1.00) -- (10.40,4.48);
  \draw[densely dashed,black!35] (12.60,1.00) -- (12.60,4.48);
  \node[font=\footnotesize,anchor=south] at (11.50,4.48) {core $U_0$};

  \draw[dotted,black!35] (8.65,3.62) -- (14.25,3.62);
  \draw[densely dashed,black!45] (8.65,2.72) -- (14.25,2.72);
  \draw[dotted,black!35] (8.65,1.82) -- (14.25,1.82);
  \node[anchor=east] at (8.42,3.62) {$\frac12+c_0\delta$};
  \node[anchor=east] at (8.42,2.72) {$\frac12$};
  \node[anchor=east] at (8.42,1.82) {$\frac12-c_0\delta$};

  \draw[blue!70!black,very thick]
    (8.72,3.62) -- (9.55,3.62)
    .. controls (9.95,3.62) and (10.00,1.82) .. (10.55,1.82)
    -- (12.45,1.82)
    .. controls (13.00,1.82) and (13.05,3.62) .. (13.45,3.62)
    -- (14.25,3.62);
  \draw[red!75!black,very thick] (8.72,3.62) -- (14.25,3.62);
  \node[text=red!75!black,anchor=south east] at (14.20,3.70)
    {$\sigma_k=+1$};
  \node[text=blue!70!black,anchor=north] at (11.50,1.75)
    {$\sigma_k=-1$};

  \draw[black!45] (8.72,0.78) -- (8.72,0.92);
  \draw[black!45] (14.25,0.78) -- (14.25,0.92);
  \node[font=\footnotesize,anchor=north] at (11.49,0.76)
    {normalized cell $U$};
  \node[font=\small] at (11.50,-0.38)
    {(b) Cell profiles};
\end{tikzpicture}
\caption{The construction of the lower bound for the statistical component. Each sign
$\sigma_k$ selects one of the two cell profiles, which agree near
$\partial Q_k$.}
\label{fig:lower-statistical-construction}
\end{figure}

For a constant $0<c_\delta\le1$ to be specified below, set
\[
\delta:=c_\delta n^{-\frac{\bar\beta}{2\bar\beta+d}},
\qquad h_j:=\delta^{1/\beta_j}.
\]
Then
\begin{equation}
V_{\bm h}=\delta^{d/\bar\beta},\qquad
n\delta^2V_{\bm h}=c_\delta^{2+d/\bar\beta}.
\label{eq:lower-statistical-scales}
\end{equation}
Define
\[
\omega_j:=\frac{\alpha\bar\beta}{d\beta_j},
\qquad
\ell_j:=h_j\left\lceil\frac{\delta^{\omega_j}}{h_j}\right\rceil.
\]
By \cref{ass:range},
\[
\sum_j\omega_j=\alpha,\qquad
0<\omega_j\le1/\beta_j,\qquad
\ell_j\asymp\delta^{\omega_j},\qquad \ell_j/h_j\in\mathbb N.
\]
Choose $c_\delta$ so that the centered box
$D_\delta=\prod_jD_{\delta,j}$ with side lengths
$\ell_j\le2\delta^{\omega_j}$ lies in $[1/4,3/4]^d$ for every $n\ge1$.
Partition it into cells $Q_1,\ldots,Q_m$ of side vector $\bm h$. Then
\begin{equation}
mV_{\bm h}
=\Leb(D_\delta)
\asymp\delta^\alpha,
\qquad m\asymp\delta^{\alpha-d/\bar\beta}\ge c>0.
\label{eq:lower-active-volume}
\end{equation}
Let $\zeta_{\pm1},U_0$ be as in \cref{lem:anisotropic-patching}, and
let $R_{D_\delta}$ be as in \cref{lem:anisotropic-barrier}.
For $\sigma\in\{-1,+1\}^m$, define
\begin{equation}
\eta_\sigma(x)
:=
\begin{cases}
\displaystyle
\frac12+c_0\delta\,
\zeta_{\sigma_k}\!\left(\operatorname{diag}(\bm h)^{-1}(x-z_k)\right),
&x\in Q_k,\\[2mm]
\displaystyle
\frac12+c_0\{\delta+R_{D_\delta}(x)\},
&x\notin D_\delta,
\end{cases}
\label{eq:clean-lower-regression}
\end{equation}
where $z_k$ is the center of $Q_k$. See
\cref{fig:lower-statistical-construction}.

Let $P_\sigma$ have regression function $\eta_\sigma$, and write
$\E_\sigma$ for expectation under $P_\sigma^{\otimes n}$.
Since $h_j^{\beta_j}=\delta$,
\cref{lem:anisotropic-patching,lem:anisotropic-barrier} give smoothness
bounds with constant $L/2$ for both the cell profiles extended by
$c_0\delta$ and the barrier $c_0R_{D_\delta}$, for sufficiently small
fixed $c_0>0$. Thus \cref{ass:holder} holds.
Taking $c_0\le1/8$ also gives $\eta_\sigma\in[1/4,3/4]$.

For $0<u<c_0\delta$, the margin event lies in $D_\delta$, so the
profile bound and \eqref{eq:lower-active-volume} give
\[
P_X\{|\eta_\sigma(X)-1/2|\le u\}
\le CmV_{\bm h}\left(\frac{u}{c_0\delta}\right)^\alpha
\le Cu^\alpha.
\]
For $c_0\delta\le u\le1$, \cref{lem:anisotropic-barrier} gives
\[
\begin{aligned}
P_X\{|\eta_\sigma(X)-1/2|\le u\}
&\le C\prod_{j=1}^d(\ell_j+u^{1/\beta_j})\\
&\le C\prod_{j=1}^d u^{\omega_j}=Cu^\alpha.
\end{aligned}
\]
Here $\ell_j\le C\delta^{\omega_j}\le Cu^{\omega_j}$ and
$u^{1/\beta_j}\le u^{\omega_j}$, since $\omega_j\le1/\beta_j$.
Thus $P_\sigma\in\cP(\alpha,\bm\beta,L)$ for sufficiently
large fixed $C_\alpha$.

For $1\le k\le m$, let $\sigma^{(k)}$ be $\sigma$ with its $k$th
sign reversed. The regressions lie in $[1/4,3/4]$ and differ only
on $Q_k$, by at most $C\delta$, so
\[
\begin{aligned}
\KL(P_\sigma^{\otimes n},P_{\sigma^{(k)}}^{\otimes n})
&\le
Cn\int_{Q_k}
\{\eta_\sigma(x)-\eta_{\sigma^{(k)}}(x)\}^2\,dx
\\
&\le
Cn\delta^2V_{\bm h}\\
&\le Cc_\delta^{2+d/\bar\beta}.
\end{aligned}
\]
For sufficiently small $c_\delta$, Pinsker's inequality gives
\begin{equation}
\TV(P_\sigma^{\otimes n},P_{\sigma^{(k)}}^{\otimes n})
\le
\sqrt{\frac12
\KL(P_\sigma^{\otimes n},P_{\sigma^{(k)}}^{\otimes n})}
\le
\frac12.
\label{eq:anisotropic-neighbor-tv}
\end{equation}
Given a classifier $\wh f$, set $C_k:=z_k+\operatorname{diag}(\bm h)U_0$ and
\[
\wh\sigma_k
:=2\1\left\{\int_{C_k}\wh f(x)\,dx
\ge\frac12\Leb(C_k)\right\}-1.
\]
Since
\[
\eta_\sigma|_{C_k}=\tfrac12+c_0\delta\sigma_k,
\qquad \Leb(C_k)=\Leb(U_0)V_{\bm h},
\]
$\wh\sigma_k\ne\sigma_k$ implies that $\wh f$ disagrees
with $f^*$ on at least half of $C_k$. Thus,
\[
\cE_0(\wh f)
\ge
c\delta V_{\bm h}
\sum_{k=1}^m\1\{\wh\sigma_k\neq\sigma_k\}
=
c\delta V_{\bm h}d_H(\wh\sigma,\sigma),
\]
where $d_H$ is the Hamming distance between $\hat\sigma$ and $\sigma$.
Applying Assouad's lemma \citep{Assouad1983,Tsybakov2009}
to each pair $(\sigma,\sigma^{(k)})$ and using
\eqref{eq:anisotropic-neighbor-tv} gives
\[
2^{-m}\sum_\sigma\E_\sigma d_H(\wh\sigma,\sigma)
\ge2^{-m}\sum_{k=1}^m2^{m-1}\cdot\frac12
=\frac m4.
\]
Using \eqref{eq:lower-decomposition},
\begin{align*}
\mathcal R_{n,r}
&\ge
\frac12\inf_{\wh f}\sup_\sigma\E_\sigma\cE_0(\wh f)\\
&\ge
c\delta V_{\bm h}\,
\inf_{\wh\sigma}2^{-m}\sum_\sigma
\E_\sigma d_H(\wh\sigma,\sigma)
\\
&\ge
c m\delta V_{\bm h}
\asymp
\delta^{1+\alpha}
\asymp
n^{-\frac{\bar\beta(1+\alpha)}{2\bar\beta+d}}.
\end{align*}

\subsection{Adversarial component}
\label{app:lower-adversarial}

When $\eta$ depends only on $t=x_1$, write
$\eta(t)$, $\Delta(t)$ and $f^*(t)$.

\begin{proof}[Proof of \cref{lem:erase-retain}]
Fix $r>0$, $z\in[0,1]^{d-1}$ and $k$.
If $f(\cdot,z)$ is constant on $J_k$, one core is misclassified. Thus
\[
\int_{J_k}
2\Delta(t)\1\{f(t,z)\neq f^*(t)\}\,dt
\ge
2c_0h_k^{\beta_-}\min\{\Leb(J_k^+),\Leb(J_k^-)\}
\ge
2c_0^2h_k^{\beta_-+1}.
\]
Otherwise, choose $b\in\partial_{J_k}\{t\in J_k:f(t,z)=1\}$.
Both values of $f(\cdot,z)$ occur arbitrarily close to $b$, so
$V_f(t,z)=1$ on
$J_k\cap(b-r/4,b+r/4)$. Since $\Leb(J_k)\ge r$,
\[
\int_{J_k} V_f(t,z)\,dt
\ge
\Leb\bigl(J_k\cap(b-r/4,b+r/4)\bigr)
\ge
\frac r4.
\]
In either case,
\[
\begin{aligned}
\int_{J_k}
\left\{
2\Delta(t)\1\{f(t,z)\neq f^*(t)\}
+
V_f(t,z)
\right\}dt
&\ge \min\left\{2c_0^2h_k^{\beta_-+1},\frac r4\right\}\\
&\ge c\min\{h_k^{\beta_-+1},r\}.
\end{aligned}
\]
Since the $J_k$ have disjoint interiors, we have
\begin{align*}
\cE_0(f)+B_r(f)
&=\int_{[0,1]^{d-1}}\int_0^1
 \left\{2\Delta(t)\1\{f(t,z)\ne f^*(t)\}+V_f(t,z)\right\}\,dt\,dz\\*
&\ge\sum_{k=1}^M\int_{[0,1]^{d-1}}\int_{J_k}
 \left\{2\Delta(t)\1\{f(t,z)\ne f^*(t)\}+V_f(t,z)\right\}\,dt\,dz\\*
&\ge c\sum_{k=1}^M\min\{h_k^{\beta_-+1},r\}.\qedhere
\end{align*}
\end{proof}

\Needspace{8\baselineskip}
We prove the adversarial bound for $0<r\le\rho$, with fixed
$\rho\in(0,1)$ chosen below. For $\rho\le r\le1$, monotonicity gives
\[
\mathcal R_{n,r}\ge\mathcal R_{n,\rho}
\ge c\psi_{\alpha,\beta_-}(\rho)
\ge c\psi_{\alpha,\beta_-}(\rho)\psi_{\alpha,\beta_-}(r),
\]
since $\psi_{\alpha,\beta_-}(r)\le1$ on $(0,1]$.
Relabel the coordinates so that $\beta_1=\beta_-$ and consider
regressions depending only on $t=x_1$.
For $M$ intervals satisfying the conditions of \cref{lem:erase-retain}
with $h_k^{\beta_-+1}\ge r$, \eqref{eq:lower-decomposition} gives
\begin{equation}
\label{eq:adversarial-block-bound}
\cE_r(f)
\ge
\frac12\{\cE_0(f)+B_r(f)\}
\ge cMr.
\end{equation}
\medskip
\noindent\emph{\bf Case 1: $\alpha\beta_-<1$.}
Set
\[
    h:=r^{1/(\beta_-+1)},
    \qquad
    M:=\left\lfloor h^{\alpha\beta_--1}\right\rfloor.
\]
Then
\[
\tfrac12h^{\alpha\beta_--1}\le M\le h^{\alpha\beta_--1},
\qquad
Mh\le h^{\alpha\beta_-}.
\]
Choose $\rho$ small enough that $6Mh\le1/2$.
Partition $D_h:=[1/4,1/4+6Mh]$ into intervals
$J_k:=[z_k-3h,z_k+3h]$. With the one-dimensional profiles and barrier of
\cref{lem:anisotropic-patching,lem:anisotropic-barrier}, define
\[
\eta(t):=\frac12+c_0
\begin{cases}
h^{\beta_-}\zeta_{-1}\!\left(\dfrac{t-z_k}{6h}\right),&t\in J_k,\\[2mm]
h^{\beta_-}+R_{D_h}(t),&t\notin D_h.
\end{cases}
\]
On each $J_k$, we use a rescaled version of the profile labelled
$\sigma_k=-1$ in \cref{fig:lower-statistical-construction}(b).
Let $P$ denote the resulting distribution.
For sufficiently small fixed $c_0$,
\cref{lem:anisotropic-patching,lem:anisotropic-barrier} ensure that
$\eta$ satisfies \cref{ass:holder} and lies in $[1/4,3/4]$,
uniformly in $h,M$.

For $0<u\le1$, \cref{lem:anisotropic-barrier}
bounds the probability of the margin event outside $D_h$ by $Cu^{1/\beta_-}$.
On each $J_k$, the profile bound gives
\[
\Leb\{t\in J_k:\Delta(t)\le u\}
\le Ch\min\left\{1,
 \left(\frac{u}{c_0h^{\beta_-}}\right)^\alpha\right\}.
\]
Summing over the blocks and adding the exterior part yields
\begin{align*}
P_X\{\Delta(X_1)\le u\}
&\le C\left[
Mh\min\left\{1,\left(\frac{u}{c_0h^{\beta_-}}\right)^\alpha\right\}
+u^{1/\beta_-}\right]\\
&\le C\{Mh^{1-\alpha\beta_-}+1\}u^\alpha\\
&\le Cu^\alpha,
\end{align*}
where $u^{1/\beta_-}\le u^\alpha$ and $Mh^{1-\alpha\beta_-}\le1$.
Thus $P\in\cP(\alpha,\bm\beta,L)$ for sufficiently large fixed
$C_\alpha$.
Each $J_k$ has the cores required by \cref{lem:erase-retain}, with
$h^{\beta_-+1}=r$ and $\Leb(J_k)=6h\ge r$.
Hence \eqref{eq:adversarial-block-bound} gives
\begin{align*}
\mathcal R_{n,r}
&\ge cMr\asymp rh^{\alpha\beta_--1}\\
&=h^{\beta_-(1+\alpha)}
=r^{\beta_-(1+\alpha)/(\beta_-+1)}.
\end{align*}
\medskip
\noindent\emph{\bf Case 2: $\alpha\beta_-=1$.}
If $\beta_-\le1$, the construction in Case 1 applies with
$h=r^{1/(\beta_-+1)}$ and $M=1$.
Thus \eqref{eq:adversarial-block-bound} gives
\[
\mathcal R_{n,r}\ge cr.
\]
For $\beta_->1$, let $t_0=1/8$ and define
\begin{equation}
g(t)
:=
t^{\beta_-}
\cos\left(2\pi\log_2\frac{t_0}{t}\right),
\qquad 0<t\le t_0,
\qquad g(0):=0.
\label{eq:log-periodic-g}
\end{equation}
For $j\ge0$, set
\[
    a_j:=2^{-j}t_0,
    \qquad
    I_j:=(a_j/2,a_j].
\]
Extend $g$ as in \cref{lem:log-periodic-model} and set
\[
    \eta(t):=\frac12+c_0g(t).
\]
The resulting model belongs to $\cP(\alpha,\bm\beta,L)$.
Each $I_j$ contains opposite-sign cores satisfying the hypotheses
of \cref{lem:erase-retain} at scale $a_j$.

Taking $\rho<t_0^{\beta_-+1}$, let
\[
M:=1+\left\lfloor
\frac{\log_2(t_0^{\beta_-+1}/r)}{\beta_-+1}
\right\rfloor.
\]
Then $a_j^{\beta_-+1}\ge r$ exactly for $0\le j<M$, and
\[
\Leb(I_j)=a_j/2\ge\tfrac12r^{1/(\beta_-+1)}\ge r.
\]
Applying \eqref{eq:adversarial-block-bound} to
$I_0,\ldots,I_{M-1}$ gives
\[
\mathcal R_{n,r}\ge crM\asymp r\log(e/r).
\]

\medskip
\noindent\emph{\bf Case 3: $\alpha\beta_->1$.}
By \cref{ass:range}, $\beta_->1$ and $\alpha\le1$.
Consider
\[
    \eta(t)
    :=
    \frac12+c_0\left(t-\frac12\right).
\]
Choose $c_0$ so that \cref{ass:holder} holds and
$\eta\in[1/4,3/4]$.
The margin condition follows from
\[
    P_X\{\Delta(X_1)\le u\}
    \le Cu
    \le Cu^\alpha,
    \qquad 0<u\le1.
\]
Set $h=1/12$ and define
\[
\begin{aligned}
J:=\left[\frac12-3h,\frac12+3h\right],
J^-:=\left[\frac12-2h,\frac12-h\right],\
\text{and}\
J^+:=\left[\frac12+h,\frac12+2h\right].
\end{aligned}
\]
Since $h\ge h^{\beta_-}$,
\[
    \Leb(J^\pm)=h,
    \qquad |\eta(t)-1/2|\ge c_0h^{\beta_-}
    \quad\text{on }J^\pm.
\]
With $\rho\le h^{\beta_-+1}$,
\eqref{eq:adversarial-block-bound} with $M=1$ gives
\[
    \mathcal R_{n,r}\ge cr.
\]
Combining the two components proves \cref{thm:lower}.

\Needspace{20\baselineskip}
\subsection{Technical lemmas for the lower bound}
\label{app:lower-technical}

\begin{lemma}
\label{lem:anisotropic-patching}
Under \cref{ass:range}, let $U=[-1/2,1/2]^d$. There exist a cube
$U_0\subset\operatorname{int}(U)$ of positive volume, functions
$\zeta_s:U\to[-1,1]$, $s\in\{-1,+1\}$, and constants $C,c_*>0$
depending only on $d,\alpha,\bm\beta,L$ such that
\[
\zeta_{+1}\equiv1,
\qquad
\zeta_{-1}=-1 \ \text{on } U_0,
\qquad
\zeta_{-1}=1 \ \text{near } \partial U,
\]
and
\[
\Leb\{u\in U:|\zeta_s(u)|\le v\}\le Cv^\alpha,
\qquad 0<v<1,\quad s\in\{-1,+1\}.
\]

Let $\bm h\in(0,1]^d$ and let
\[
Q_k=z_k+\operatorname{diag}(\bm h)U\subset[0,1]^d,
\qquad k=1,\ldots,N,
\]
have disjoint interiors, and let $s_k\in\{-1,+1\}$. If
$
0<a\le c_*\min_{1\le j\le d} h_j^{\beta_j},
$
then
\[
F(x):=
\begin{cases}
a\zeta_{s_k}\!\left(
\operatorname{diag}(\bm h)^{-1}(x-z_k)
\right),
& x\in Q_k,\\[2mm]
a,
& x\notin\displaystyle\bigcup_{k=1}^N Q_k,
\end{cases}
\]
satisfies \cref{ass:holder}, uniformly in
$N$, $\bm h$, and $(s_1,\ldots,s_N)$.
\end{lemma}

\begin{proof}
Choose a smooth radial $q_0:U\to[-1,1]$, equal to $-1$ near the
origin and $1$ for $\|u\|_2\ge1/4$, with a single simple radial zero.
Take $U_0\subset\{q_0=-1\}$. Then
\[
\Leb\{u\in U:|q_0(u)|\le w\}\le Cw,
\qquad 0<w<1.
\]
Set
\[
\zeta_{+1}\equiv1,
\qquad
\zeta_{-1}:=
\begin{cases}
q_0,
& \alpha\le1,\\[1mm]
\sgn(q_0)|q_0|^{1/\alpha},
& \alpha>1.
\end{cases}
\]

If $\alpha>1$, \cref{ass:range} implies
$\beta_j\le1/\alpha<1$ for every $j$, and $\zeta_{-1}$ is
$1/\alpha$-H\"older.
The sublevel bound follows from that for $q_0$ and the
definition of $\zeta_{-1}$. For $\zeta_{+1}$, the sublevel
set is empty for $0<v<1$.

Extend $\zeta_s-1$ by zero outside $U$. Coordinate derivatives
match across cell faces, since it vanishes near $\partial U$.
Let $F_\nu$, $1\le\nu\le N$, be the zero extension from $Q_\nu$ of
$a[\zeta_{s_\nu}\!\left(\operatorname{diag}(\bm h)^{-1}(x-z_\nu)\right)-1]$.
If $\alpha\le1$, smoothness of the profiles gives
\[
\begin{aligned}
\|\partial_j^kF_\nu\|_\infty
&\le Ca h_j^{-k}\le Cc_*,\qquad 1\le k\le m_j,\\
|\partial_j^{m_j}F_\nu(x+te_j)-\partial_j^{m_j}F_\nu(x)|
&\le Ca h_j^{-m_j}\min\{1,|t|/h_j\}\\
&\le Cc_*|t|^{\beta_j-m_j}.
\end{aligned}
\]
Here we used $a\le c_*h_j^{\beta_j}$.
If $\alpha>1$, then $m_j=0$ and $\beta_j\le1/\alpha$.
The profile then gives
\[
|F_\nu(x+te_j)-F_\nu(x)|
\le Ca\min\{1,(|t|/h_j)^{1/\alpha}\}
\le Cc_*|t|^{\beta_j}.
\]
In $F=a+\sum_\nu F_\nu$, at most one $F_\nu$ contributes at a point
and at most two to an increment. Choosing $c_*$ small enough gives
\cref{ass:holder}, uniformly in $N$, $\bm h$ and the signs.
\end{proof}

\Needspace{12\baselineskip}
\begin{lemma}
\label{lem:anisotropic-barrier}
For every box $D=\prod_jD_j\subset(0,1)^d$, there exists
$R_D:[0,1]^d\to[0,1]$ that vanishes on $D$ and satisfies the
coordinate derivative and H\"older bounds in \cref{ass:holder}
with a constant $C$ depending only on $d,\bm\beta$.
Moreover,
\[
\dist(x_j,D_j)\le C R_D(x)^{1/\beta_j},
\qquad x\in[0,1]^d,\quad j=1,\ldots,d.
\]
\end{lemma}

\Needspace{5\baselineskip}
\begin{proof}
Set
\[
R_D(x):=\frac1d\sum_{j=1}^d\dist(x_j,D_j)^{\beta_j}.
\]
Since each coordinate distance is at most $1$, we have
$0\le R_D\le1$, and $R_D=0$ on $D$.
The distance bound follows from
$\dist(x_j,D_j)^{\beta_j}\le dR_D(x)$.

Write $D_j=[a_j,b_j]$. Then
\[
\dist(t,D_j)^{\beta_j}
=(a_j-t)_+^{\beta_j}+(t-b_j)_+^{\beta_j},
\qquad t\in[0,1],
\]
where $u_+:=\max\{u,0\}$.
Recall that $m_j=\lceil\beta_j\rceil-1<\beta_j$.
Both one-sided powers have $m_j$ continuous derivatives, and their
$m_j$th derivatives are constant multiples of one-sided powers
of order $\beta_j-m_j\in(0,1]$.
These derivatives satisfy the required H\"older bound uniformly in $D_j$.
All derivatives up to order $m_j$ are also uniformly bounded on
$[0,1]$, with constants depending only on $\beta_j$.
Along coordinate $j$, only the $j$th summand in $R_D$ varies,
which gives the required regularity bounds.
\end{proof}

\Needspace{12\baselineskip}
\begin{lemma}
\label{lem:log-periodic-model}
Suppose that $\beta_->1$ and $\alpha=1/\beta_-$.
There exists an extension of $g$ in \eqref{eq:log-periodic-g} to $[0,1]$
such that the uniform design and regression
$\eta(t)=1/2+c_0g(t)$, with $t=x_1$, define a distribution in
$\cP(\alpha,\bm\beta,L)$ for sufficiently small fixed $c_0>0$
and sufficiently large $C_\alpha$.

For the intervals $I_j$ defined in \cref{app:lower-adversarial},
there exist disjoint subintervals $I_j^+,I_j^-\subset I_j$ such that
\[
\Leb(I_j^+)\wedge\Leb(I_j^-)\ge ca_j,
\qquad
g\ge ca_j^{\beta_-}\ \text{on }I_j^+,
\qquad
g\le-ca_j^{\beta_-}\ \text{on }I_j^-.
\]
\end{lemma}

\begin{proof}
Choose a smooth cutoff $0\le\chi\le1$, equal to one on $[0,t_0]$
and zero on $[9t_0/8,1]$. Extend $g$ by
\[
g(t):=\chi(t)t^{\beta_-}
\cos\left(2\pi\log_2\frac{t_0}{t}\right)
+\{1-\chi(t)\}t_0^{\beta_-},
\qquad 0<t\le1.
\]
This extension is smooth for $t>0$ and bounded away from zero
on $[t_0,1]$.

Let $m:=\lceil\beta_-\rceil-1$ and $\vartheta:=\beta_--m\in(0,1]$.
For $0<t\le t_0$, $g^{(k)}(t)$ is $t^{\beta_--k}$ times a smooth
bounded periodic function of $\log t$. Consequently,
$g^{(k)}(0)=0$ for $k\le m$, and
\[
|g^{(m)}(t)|\le Ct^\vartheta,
\qquad |g^{(m+1)}(t)|\le Ct^{\vartheta-1}.
\]
For $0\le s<t\le t_0$, these bounds give
\[
|g^{(m)}(t)-g^{(m)}(s)|\le
\begin{cases}
C(t^\vartheta+s^\vartheta)\le C(t-s)^\vartheta,&s\le t/2,\\[1mm]
Ct^{\vartheta-1}(t-s)\le C(t-s)^\vartheta,&s>t/2.
\end{cases}
\]
Thus $g$ satisfies the H\"older condition of order $\beta_-$.

On $I_j$, substitute $t=a_j2^{-v}$, $0\le v<1$. Then
\[
g(a_j2^{-v})=a_j^{\beta_-}2^{-\beta_-v}\cos(2\pi v),
\qquad
\left|\frac{dt}{dv}\right|=(\log2)a_j2^{-v}\asymp a_j.
\]
Since $\cos(2\pi v)$ has simple zeros, this change of
variables gives
\[
\Leb\{t\in I_j:|g(t)|\le u\}
\le Ca_j\min\{1,u/a_j^{\beta_-}\}.
\]
For $0<u\le1$, dyadic summation and $\beta_->1$ give
\[
\begin{aligned}
\sum_{j\ge0}a_j\min\{1,u/a_j^{\beta_-}\}
&=\sum_{j:a_j\le u^{1/\beta_-}}a_j
 +u\sum_{j:a_j>u^{1/\beta_-}}a_j^{1-\beta_-}\\
&\le Cu^{1/\beta_-}.
\end{aligned}
\]
Together with the positive lower bound for $g$ on
$[t_0,1]$, this yields
\[
\Leb\{t\in[0,1]:|g(t)|\le u\}\le Cu^{1/\beta_-}.
\]
Choose $c_0$ small enough for \cref{ass:holder} and $0\le\eta\le1$.

Take
\[
I_j^+:=[a_j2^{-1/6},a_j],
\qquad I_j^-:=[a_j2^{-2/3},a_j2^{-1/3}].
\]
Both intervals have length comparable to $a_j$. The cosine factor
is at least $1/2$ on $I_j^+$ and at most $-1/2$ on $I_j^-$.
\end{proof}

\begin{proof}[Proof of \cref{cor:population}]
For each fixed distribution, any measurable classifier has expected
excess risk at least $\inf_f\cE_r(f)$. Hence \cref{cor:minimax} gives
\[
\sup_{P\in\cP}\inf_f\cE_r(f)
\le\mathcal R_{n,r}
\le C\left\{
n^{-\frac{\bar\beta(1+\alpha)}{2\bar\beta+d}}
+\psi_{\alpha,\beta_-}(r)\right\}.
\]
Letting $n\to\infty$ with $r>0$ fixed proves the upper bound.
For the lower bound, the constructions above and
\eqref{eq:adversarial-block-bound} give, for every sufficiently small
$r>0$, a distribution in $\cP$ such that
$\cE_r(f)\ge c\psi_{\alpha,\beta_-}(r)$ for every classifier $f$.
Taking the infimum over $f$ and the supremum over the class proves
the lower bound. The extension to the remaining radii follows by
monotonicity, as above.
\end{proof}

\section{Proofs for adaptive classification}
\label{app:adaptive}

\subsection{A holdout oracle inequality}
\label{app:holdout-oracle}

\begin{proof}[Proof of \cref{lem:bernstein}]
Since both losses are binary, $|W_{r,f}(X,Y)|\le1$. If
$V_f(X)=0$ and $f(X)=f^*(X)$, the two losses agree. Hence
\[
W_{r,f}(X,Y)^2
\le V_f(X)+\1\{f(X)\ne f^*(X)\},
\]
so
\[
\operatorname{Var}\{W_{r,f}(X,Y)\}
\le B_r(f)+P_X\{f(X)\ne f^*(X)\}.
\]
For every $u>0$, the margin condition gives
\begin{align*}
P_X\{f(X)\ne f^*(X)\}
&\le P_X\{\Delta(X)\le u\}
 +\frac1u\E\{\Delta(X)\1\{f(X)\ne f^*(X)\}\}\\
&\le C_\alpha u^\alpha+\frac{\cE_0(f)}{2u}.
\end{align*}
By \cref{thm:decomposition},
$B_r(f)\le2\cE_r(f)$ and $\cE_0(f)\le2\cE_r(f)$. Then, we have
\[
\operatorname{Var}\{W_{r,f}(X,Y)\}
\le2\cE_r(f)+C_\alpha u^\alpha+\frac{\cE_r(f)}{u}.
\]
If $\cE_r(f)>0$, take $u=\cE_r(f)^{1/(1+\alpha)}$ and use $\cE_r(f)\le1$ to obtain
\[
\operatorname{Var}\{W_{r,f}(X,Y)\}
\le(C_\alpha+3)\cE_r(f)^{\alpha/(1+\alpha)}.
\]
If $\cE_r(f)=0$, the variance is at most $C_\alpha u^\alpha$ for every
$u>0$, so letting $u\rightarrow0$ proves the claim.
\end{proof}

\begin{proof}[Proof of \cref{thm:holdout}]
Conditional on $\cD_1$, the candidates are fixed and the validation
observations remain independent with common distribution $P$.
Let $P_2$ denote the empirical measure of $\cD_2$.
If $\log(2N/\delta)\le n_2$, with probability
at least $1-\delta$, \cref{lem:bernstein} and Bernstein's inequality give
\begin{equation}
\left|(P_2-P)W_{r,\wh f_m}\right|
\le\tau\cE_r(\wh f_m)
+C_\tau\left\{\frac{\log(2N/\delta)}{n_2}\right\}^{\frac{1+\alpha}{2+\alpha}},
\qquad 1\le m\le N.
\label{eq:holdout-deviation}
\end{equation}
When $\log(2N/\delta)>n_2$, the probability bound follows
directly from $\cE_r(\wh f_r)\le1$ after increasing $C_\tau$.
Hence it holds for every $\delta\in(0,1)$.

Since $\1\{f^*(x)\ne y\}$ does not depend on $m$, the holdout
selector also minimizes $P_2W_{r,\wh f_m}$ over $m$.
Applying \eqref{eq:holdout-deviation} to the selected classifier
and each candidate gives
\[
(1-\tau)\cE_r(\wh f_r)
\le(1+\tau)\cE_r(\wh f_m)
+C_\tau\left\{\frac{\log(2N/\delta)}{n_2}\right\}^{\frac{1+\alpha}{2+\alpha}}.
\]
Taking the minimum over $m$ and dividing by $1-\tau$ proves
the probability bound.

For the expectation bound, take $\tau=1/3$ and
$\delta=e^{-t}$ for $t>0$. Integrating the conditional tail bound gives
\begin{align*}
\E\{\cE_r(\wh f_r)\mid\cD_1\}
&\le2\min_{1\le m\le N}\cE_r(\wh f_m)+C\int_0^\infty
\left\{\frac{\log(2N)+t}{n_2}\right\}^{\frac{1+\alpha}{2+\alpha}}
e^{-t}\,dt\\
&\le2\min_{1\le m\le N}\cE_r(\wh f_m)
+C\left\{\frac{\log(2N)}{n_2}\right\}^{\frac{1+\alpha}{2+\alpha}}.
\end{align*}
Taking expectation over $\cD_1$ completes the proof.
\end{proof}

\subsection{Adaptive bandwidth selection}
\label{app:bandwidth-adaptation}

\begin{proof}[Proof of \cref{thm:adaptive}]
For $0\le r\le1$, let $h_j$ be the smallest dyadic
bandwidth satisfying
\[
h_j\ge
n_1^{-\frac{\bar\beta}{\beta_j(2\bar\beta+d)}}
\vee r^{\frac{\beta_-}{\beta_j(\beta_-+1)}},
\qquad j=1,\ldots,d.
\]
In particular,
\[
\min_jh_j\ge r^{1/(\beta_-+1)}\ge r,
\qquad
n_1V_{\bm h}
\ge n_1^{\frac{2\bar\beta}{2\bar\beta+d}}\ge1.
\]
For $\bm h=2^{-\bm k}$, the bound $n_1V_{\bm h}\ge1$ implies
$2^{|\bm k|}=V_{\bm h}^{-1}\le n_1$, and hence $|\bm k|\le J_n$.
Since $\min_jh_j\ge r$, we have $\bm h\in\mathcal H_n(r)$.
The number of candidates satisfies
\[
|\mathcal H_n(r)|\le|\mathcal H_n|
=\binom{J_n+d}{d}\le C(\log n)^d.
\]
Using \cref{thm:holdout,thm:fixed-upper}, we obtain
\begin{align*}
\sup_{P\in\cP(\alpha,\bm\beta,L)}
\E\cE_r(\wh f_{\mathrm{ad},r})
&\le2\sup_{P\in\cP(\alpha,\bm\beta,L)}
\E\cE_r(\wh f_{\bm h}^{(1)})
+C\left(\frac{\log n\,\log\log n}{n}\right)^{\frac{1+\alpha}{2+\alpha}}\\
&\le C\left\{
n^{-\frac{\bar\beta(1+\alpha)}{2\bar\beta+d}}
+\psi_{\alpha,\beta_-}(r)
+\left(\frac{\log n\,\log\log n}{n}\right)^{\frac{1+\alpha}{2+\alpha}}
\right\}.
\end{align*}
If $\alpha\bar\beta<d$, then
\[
\frac{1+\alpha}{2+\alpha}
>\frac{\bar\beta(1+\alpha)}{2\bar\beta+d}.
\]
Consequently,
\[
\left(\frac{\log n\,\log\log n}{n}\right)^{\frac{1+\alpha}{2+\alpha}}
\le C n^{-\frac{\bar\beta(1+\alpha)}{2\bar\beta+d}}.
\]
The selection remainder is therefore absorbed into the statistical
term with a constant independent of $r$.
\end{proof}

\section{Proofs under an ordinary-risk constraint}
\label{app:budget-proof}

\begin{proof}[Proof of \cref{thm:budget-frontier}]
Let $\delta:=n^{-\bar\beta/(2\bar\beta+d)}$.
All bounds below are uniform over $P\in\cP$.

\emph{Upper bound.}
For $\underline h^{-1}\in\mathbb N$ and
$\delta^{1/\beta_-}\le\underline h\le1$, take
$h_j=1/\lfloor\underline h^{-\beta_-/\beta_j}\rfloor$.
Then, $\min_jh_j=\underline h$, and the bounds in
\cref{app:upper-completion} give
\[
b_{\bm h}+v_{\bm h}\le C\underline h^{\beta_-},\qquad
nV_{\bm h}\ge\underline h^{-2\beta_-}\ge1.
\]
Thus \eqref{eq:ordinary-final} yields
\begin{equation}
\sup_{P\in\cP}\E\cE_0(\wh f_{\bm h})
\le C\underline h^{\beta_-(1+\alpha)}.
\label{eq:budget-ordinary}
\end{equation}
We also have
\begin{equation}
\sup_{P\in\cP}\E B_r(\wh f_{\bm h})
\le C\min\{r^{\alpha\beta_-},\,r\underline h^{\alpha\beta_--1}\}.
\label{eq:budget-exposure}
\end{equation}
For $r\le\underline h$, \eqref{eq:rough-G-sharp} gives
$Cr\underline h^{\alpha\beta_--1}\le Cr^{\alpha\beta_-}$.
For $r>\underline h$ and $\Delta(x)>C_*r^{\beta_-}$,
\eqref{eq:coordinate-modulus} implies that $f^*$ is constant
on $A_r(x)$ and that $\Delta(x')\ge\Delta(x)/2$ for every
$x'\in A_r(x)$.
An attack can change the prediction at $x$ only if
$\|\wh P_Q-\eta\|_{\infty,Q}\ge\Delta(x)/2$ on one of at most
$C(r/\underline h)^d$ cells meeting $A_r(x)$.
Choose $C_*$ so that $C_{\rm app}b_{\bm h}\le\Delta(x)/4$.
Applying \cref{lem:cell-concentration} with $t=\Delta(x)/8$ gives
\[
P^{\otimes n}\{V_{\wh f_{\bm h}}(x)=1\}
\le C(r/\underline h)^d e^{-c\Delta(x)^2/\underline h^{2\beta_-}}.
\]
Integration over the dyadic margin layers above $C_*r^{\beta_-}$ yields
\begin{align*}
\E B_r(\wh f_{\bm h})
&\le Cr^{\alpha\beta_-}
 \left\{1+(r/\underline h)^d
 \sum_{k\ge0}2^{k\alpha}e^{-c4^k(r/\underline h)^{2\beta_-}}\right\}\\
&\le Cr^{\alpha\beta_-}.
\end{align*}
For a sufficiently large fixed $C_1\ge1$, choose the smallest
$\underline h$ with $\underline h^{-1}\in\mathbb N$ satisfying
\[
\underline h\ge\delta^{1/\beta_-}\vee\min\left\{
 (\varepsilon/(2C_1))^{1/[\beta_-(1+\alpha)]},\,
 r^{1/(\beta_-+1)}\right\}.
\]
Taking $C_0\ge2C_1$, \eqref{eq:budget-ordinary} gives
\begin{equation}
\sup_{P\in\cP}\E\cE_0(\wh f_{\bm h})
\le C\underline h^{\beta_-(1+\alpha)}
\le\max\{C_1\delta^{1+\alpha},\varepsilon/2\}
\le\varepsilon/2.
\label{eq:budget-feasibility}
\end{equation}
By \cref{thm:decomposition}, it remains to evaluate the exposure term.
If $\varepsilon\le r^{\beta_-(1+\alpha)}$, then
\[
\E\cE_r(\wh f_{\bm h})
\le\varepsilon/2+Cr^{\alpha\beta_-}\le Cr^{\alpha\beta_-}.
\]
If $r^{\beta_-(1+\alpha)}<\varepsilon\le r^{\beta_-(1+\alpha)/(\beta_-+1)}$, then
$\underline h\asymp\varepsilon^{1/[\beta_-(1+\alpha)]}$, so
\[
\E\cE_r(\wh f_{\bm h})
\le\varepsilon/2+Cr\varepsilon^{-\frac{1-\alpha\beta_-}{\beta_-(1+\alpha)}}
\le Cr\varepsilon^{-\frac{1-\alpha\beta_-}{\beta_-(1+\alpha)}},
\]
where $\varepsilon^{\frac{\beta_-+1}{\beta_-(1+\alpha)}}\le r$.
If $\varepsilon>r^{\beta_-(1+\alpha)/(\beta_-+1)}$, then
$\underline h\asymp\delta^{1/\beta_-}\vee r^{1/(\beta_-+1)}$, and
\eqref{eq:budget-ordinary}--\eqref{eq:budget-exposure} give
\[
\E\cE_r(\wh f_{\bm h})
\le C\{\delta^{1+\alpha}+r^{\beta_-(1+\alpha)/(\beta_-+1)}\}.
\]

\Needspace{5\baselineskip}
\emph{Lower bound.}
\Cref{thm:lower} supplies the statistical term in every regime
and the entire lower bound for $\varepsilon>r^{\beta_-(1+\alpha)/(\beta_-+1)}$.
Suppose $r>0$ and $\varepsilon\le r^{\beta_-(1+\alpha)/(\beta_-+1)}$.
We first use the construction in \cref{app:lower-adversarial}, with
\[
h:=K\max\{r,\varepsilon^{\frac{1}{\beta_-(1+\alpha)}}\},\qquad
M:=\lfloor h^{\alpha\beta_--1}\rfloor.
\]
On each block, a constant prediction contributes at least
$c_1h^{\beta_-+1}$ to ordinary excess risk.
Choose $K\ge1$ so that $c_1K^{\beta_-(1+\alpha)}\ge4$.
Since $h\le Kr^{1/(\beta_-+1)}$, there is a fixed $\rho\in(0,1)$
such that the packing conditions hold for $0<r\le\rho$, with
$r\le h$ and $M\ge\tfrac12h^{\alpha\beta_--1}$.
In particular,
\[
c_1Mh^{\beta_-+1}\ge(c_1/2)h^{\beta_-(1+\alpha)}\ge2\varepsilon.
\]
On each block and coordinate slice, the argument in the proof of
\cref{lem:erase-retain} gives either ordinary excess loss at least
$c_1h^{\beta_-+1}$ or exposed length at least $r/4$.
Summing the corresponding inequalities over the blocks and integrating over the slices gives
\[
B_r(f)\ge\frac r4
 \left\{M-\frac{\cE_0(f)}{c_1h^{\beta_-+1}}\right\}.
\]
Then, every classifier satisfying the ordinary-risk budget has
\[
\E B_r(\wh f)\ge\frac r4
 \left\{M-\frac{\varepsilon}{c_1h^{\beta_-+1}}\right\}
\ge\frac{rM}{8}
\ge c\min\{r^{\alpha\beta_-},\,r\varepsilon^{-\frac{1-\alpha\beta_-}{\beta_-(1+\alpha)}}\}.
\]
Combining this with \cref{thm:decomposition,thm:lower} proves the
lower bound for $0\le r\le\rho$.
For $\rho<r\le1$, monotonicity in $r$ and \cref{thm:lower} give
\[
\mathcal R_{n,r}(\varepsilon)
\ge\mathcal R_{n,r}
\ge\mathcal R_{n,\rho}
\ge c\rho^{\beta_-(1+\alpha)/(\beta_-+1)}.
\]
The right-hand side of \eqref{eq:budget-frontier} is at most $2$,
so adjusting the comparison constant proves the lower bound over
the remaining radii as well.
\end{proof}

\begin{proof}[Proof of \cref{cor:budget-conflict}]
Let $\theta:=\beta_-(1+\alpha)/(\beta_-+1)$.
The assumptions imply
$\alpha\beta_-<\theta<1$ and $\varepsilon_n=o(r_n^\theta)$.
Hence \cref{cor:minimax} gives
$\mathcal R_{n,r_n}\asymp r_n^\theta$.
For all sufficiently large $n$, the budget lies in the tight
or moderate regime of \cref{thm:budget-frontier}. Therefore, since $r_n\to0$ and $\varepsilon_n=o(r_n^\theta)$,
\[
\frac{\mathcal R_{n,r_n}(\varepsilon_n)}
     {\mathcal R_{n,r_n}}
\gtrsim
\min\left\{
r_n^{\alpha\beta_--\theta},\,
\left(\frac{r_n^\theta}{\varepsilon_n}\right)^{
   (1-\theta)/\theta}
\right\}
\longrightarrow\infty.
\]
\end{proof}

\bibliographystyle{plainnat}
\bibliography{ref}

\end{document}